\documentclass[dvipsnames]{article}

\usepackage[normalem]{ulem}
\usepackage{a4wide}

\usepackage{amsmath,amsfonts,amsthm,amssymb,bbm}
\usepackage[justification=justified, singlelinecheck=false]{caption}
\usepackage{graphicx}
\usepackage{enumitem}
\usepackage{setspace}

\usepackage{tikz}

\usepackage{algorithmicx}
\usepackage{algorithm}
\usepackage{algpseudocode}

\allowdisplaybreaks

\usepackage[scaled=.90]{helvet}
\usepackage{courier}
\usepackage{ae}
\usepackage[T1]{fontenc}
\usepackage{subfigure}
\usepackage{graphicx,color}
\usepackage{xcolor}

\usepackage[T1]{fontenc}
\usepackage{here} 
\usepackage[colorlinks=false, pdfborder={0 0 0}]{hyperref}

\usepackage{xpatch}
\makeatletter
\xpatchcmd{\@thm}{\thm@headpunct{.}}{\thm@headpunct{}}{}{}
\makeatother

\newcommand{\norm}[1]{\left\lVert#1\right\rVert}

\newtheorem{theorem}{Theorem}

\newtheorem{lemma}{Lemma}

\newtheorem{remark}{Remark}
\newtheorem{proposition}{Proposition}

\theoremstyle{definition}

\begin{document}

\title{Wiener-path-based exact and splitting simulation\\ of linear and semi-linear SDEs with additive noise}

\author{
Irene Tubikanec \\
Institute of Applied Statistics, Johannes Kepler University Linz, Austria \\[-0.4em]
\texttt{irene.tubikanec@jku.at}
}
\date{}
\maketitle	
	
\thispagestyle{empty}

\vspace{0.4cm}

\section*{Abstract}

Strong, or pathwise, simulation of stochastic differential equations (SDEs) requires linking a generated SDE path to the corresponding realization of its driving Wiener process. In this work, we develop Wiener-path-based simulation procedures for multidimensional linear and semi-linear SDEs with additive noise that preserve this connection.

Building upon a known Wiener-path-based exact simulation approach for linear SDEs, we extend this construction to generate coherent paths simultaneously on different time grids. We further propose an alternative exact simulation procedure that requires fewer random numbers and therefore provides a more efficient construction, while retaining a controlled approximation of the coupling between the SDE path and its associated Wiener path. The proposed methods are illustrated on a stochastic harmonic oscillator, a well-known model from classical mechanics, with the simulation results confirming the theoretical findings.

We then extend the Wiener-path-based construction to semi-linear SDEs, which arise in many real-world applications. In contrast to linear SDEs, their solutions are generally not available in closed form, making numerical approximation methods necessary. For Euler--Maruyama-type methods, the relation between a given Wiener path and the corresponding approximate SDE path is explicit, whereas constructing Wiener-path-based simulations using structure-preserving splitting methods is less immediate. We address this challenge by incorporating the aforementioned Wiener-path-based constructions into splitting methods. Illustrations on the stochastic FitzHugh--Nagumo model and the stochastic Jansen--Rit neural mass model, two important models from neuroscience, demonstrate the applicability of the proposed procedures to non-linear stochastic systems.

\subsubsection*{Keywords} Strong pathwise simulation, Linear SDEs, Exact simulation, Semi-linear SDEs, Numerical splitting methods. 

\subsubsection*{AMS subject classifications} 60H10, 60H35, 65C30


\vspace{0.5cm}


\section{Introduction}
\label{sec1:Intro}

This manuscript considers $d$-dimensional semi-linear SDEs with additive noise of the form 
\begin{equation*}
    dX(t)=\left[ AX(t)+N(X(t)) \right] dt+\Sigma dW(t), \quad t  \in[0,T], \quad X(0)=X_0,
\end{equation*}
where $A\in\mathbb{R}^{d\times d}$, $\Sigma\in\mathbb{R}^{d\times m}$, $(W(t))_{t\in [0,T]}$ is an $m$-dimensional standard Wiener process, and the drift function $N\colon \mathbb{R}^d \to \mathbb{R}^d$ is either globally Lipschitz continuous or satisfies a global one-sided Lipschitz condition while being allowed to exhibit polynomial growth.

Let $0=t_0<t_1<\ldots<t_{n-1}<t_n=T$ with $t_i=i\Delta$, $i=0,\ldots,n$, and equidistant step-size $\Delta$. For the standard Euler--Maruyama method \cite{Kloeden1992}, as well as its tamed \cite{Hutzenthaler2012_2,Sabanis2013,Tretyakov2012,Zhang2017} and truncated \cite{Mao2017,Hutzenthaler2012,Mao2015,Mao2016} variants required in the non-globally Lipschitz setting, the relation between a given Wiener path $W(t_i)$, $i=0,\ldots,n$, and the corresponding approximate SDE path $\widetilde{X}(t_i)\approx X(t_i)$, $i=0,\ldots,n$, is explicit. This is because these methods directly involve Wiener increments
\begin{equation*}
   \xi_i(\Delta) := W(t_i)-W(t_{i-1})=\int\limits_{t_{i-1}}^{t_i} \mathbb{I}_{m\times m} \ dW(s), \qquad i=1,\ldots,n,
\end{equation*}
which are independent and have a multivariate normal distribution with zero mean vector and covariance matrix 
\begin{equation*}
    C_\xi(\Delta) :=\textrm{Cov}\left( \xi_i(\Delta) \right) = \Delta \mathbb{I}_{m\times m},
\end{equation*}
where $\mathbb{I}_{m\times m}$ denotes the $m\times m$-dimensional identity matrix.

However, it is well-known that Euler--Maruyama-type methods may fail to preserve important structural properties of SDEs; see, e.g., \cite{Ableidinger2017,Buckwar2022,Chevallier2020,Cohen2012,Ditlevsen2019,Mattingly2002} and the references therein. In contrast, splitting methods are often specifically designed to preserve qualitative features of the underlying dynamics. We refer to \cite{Blanes2009,Mclachlan2002} for splitting methods of ordinary differential equations (ODEs) and to \cite{Ableidinger2016,Alamo2016,Brehier2018,Chevallier2020,Fosteretal2022,JansLord2025,KellyLord2023,Misawa2001,Petersen1998,Shardlow2003,Skeel1999,Tubikanec2020} for selected extensions to SDEs. Specifically, splitting methods for semi-linear SDEs as above are, e.g., discussed in \cite{Ableidinger2017,Buckwar2022,Ditlevsenetal2025,Pilipovic2024}.

Unlike for Euler--Maruyama-type methods, it is less immediate how to associate a given Wiener path with the corresponding approximate SDE path when using splitting methods. The reason is that such methods, for semi-linear SDEs, typically involve stochastic integrals of the form
\begin{equation*}
    \eta_i(\Delta):=\int\limits_{t_{i-1}}^{t_i} e^{A(t_i-s)} \Sigma \ dW(s), \qquad i=1,\ldots,n,
\end{equation*}
which are independent multivariate normal random variables with zero mean vector and covariance matrix 
\begin{equation*}
    C_\eta(\Delta):=\textrm{Cov}\left( \eta_i(\Delta) \right)=\int\limits_{t_{i-1}}^{t_i} e^{A(t_i-s)} \Sigma\Sigma^\top \left( e^{A(t_i-s)} \right)^\top \ ds. 
\end{equation*}
Consequently, standard simulation procedures for splitting methods generate realizations of $\eta_i(\Delta)$ directly from the corresponding multivariate normal distribution, e.g., using a Cholesky factorization of $C_\eta(\Delta)$. However, this does not provide an approximate SDE path that is coherent with a specific underlying Wiener path. 

The aim of this manuscript is to address this limitation by proposing a simulation procedure that generates approximate paths of the SDE via splitting methods, which are coherent with their specific underlying Wiener paths. To achieve this, we build upon the approach proposed in \cite{Cruz2020} for linear SDEs with additive noise of the form
\begin{equation*}
    dX(t)=AX(t) dt+\Sigma dW(t), \quad t  \in[0,T], \quad X(0)=X_0,
\end{equation*}
where $A$, $\Sigma$, and $(W(t))_{t\in [0,T]}$ are as above. The unique strong solution $(X(t))_{t\in [0,T]}$ to this equation is known exactly and, considered at the grid points, also involves the stochastic integrals $\eta_i(\Delta)$, $i=1,\ldots,n$. Instead of using only their covariance matrix $C_\eta(\Delta)$, the authors derive the joint covariance matrix of $\xi_i(\Delta)$ and $\eta_i(\Delta)$, denoted here by $C_{\xi;\eta}(\Delta)$, and jointly sample from $(\xi_i(\Delta),\eta_i(\Delta))^\top$, for $i=1,\ldots,n$, to construct a path $X(t_i)$, $i=0,\ldots,n$, of the linear SDE that is linked to its underlying Wiener path $W(t_i)$, $i=0,\ldots,n$. In~this manuscript, we recall this \textit{Wiener-path-based exact simulation} procedure for linear SDEs and combine it with existing splitting methods for semi-linear SDEs, resulting in the proposed \textit{Wiener-path-based splitting simulation} algorithm. 

Another aspect addressed in this work is the simultaneous generation of an SDE path on two different time grids. For Euler--Maruyama-type methods, this is straightforward and can be achieved by subsampling the given Wiener path on the fine grid to obtain the corresponding path on the coarse grid. For splitting methods, however, such a construction is less immediate. Building again upon \cite{Cruz2020}, we therefore first propose a \textit{Wiener-path-based exact simulation on different time grids} procedure for generating a path of a linear SDE simultaneously on a fine and a coarse time grid. We then extend this approach to \textit{Wiener-path-based splitting simulation on different time grids} for semi-linear SDEs. 

Finally, this work investigates an alternative simulation approach for the case in which the SDE and its driving Wiener process have the same dimension, i.e., $d=m$. The proposed approach requires fewer random numbers and is therefore more efficient than the aforementioned Wiener-path-based construction. Specifically, the latter requires joint sampling from $(\xi_i(\Delta),\eta_i(\Delta))^\top$, which amounts to generating realizations of an $(m+d)$-dimensional normal distribution. In contrast, for $d=m$, our alternative simulation procedure requires sampling only from a single $d$-dimensional normal distribution. This construction preserves the correct marginal distributions of $\xi_i(\Delta)$ and $\eta_i(\Delta)$, but approximates their cross-covariance structure. We show that the resulting SDE path remains linked to its underlying Wiener path through an approximate cross-covariance structure, whose deviation from the exact cross-covariance is of order $\Delta^2$ for $d=m=1$ and of order $\Delta$ for $d=m>1$. This alternative approach is first introduced and theoretically investigated for linear SDEs and then adapted to splitting simulation for semi-linear~SDEs.

The contributions of this manuscript are summarized as follows:
\begin{itemize}
    \item In Section \ref{sec2:LinearSDEs}, we recall basic results on linear SDEs and, following \cite{Cruz2020}, review a construction that enables the generation of a solution path coherent with a specific underlying Wiener path. We then extend this construction to enable the generation of coherent paths on different time grids. 
    \item In Section \ref{sec3:simple_SHO}, we illustrate these procedures on a  stochastic harmonic oscillator, a well-known linear SDE arising in classical mechanics. 
    \item In Section \ref{sec4:SemiLinearSDEs}, we show how these results can be incorporated into splitting methods for semi-linear SDEs and present the resulting Wiener-path-based splitting simulation algorithms. 
    \item In Section~\ref{sec5:Neuronal_models}, we illustrate the proposed algorithms on two stochastic models from neuroscience of practical relevance: the stochastic Fitz-Hugh Nagumo (FHN) model and the stochastic Jansen--Rit neural mass model (JRNMM).
    \item In Section \ref{appA}, we present the alternative simulation procedure for the case $d=m$. These results are first discussed for linear SDEs and illustrated on the Ornstein--Uhlenbeck process and the stochastic harmonic oscillator. We then adapt the approach to splitting simulation for non-linear SDEs and illustrate it on the stochastic FHN model.
\end{itemize}
The manuscript concludes with Section \ref{sec:conclusion}, which summarizes the main findings. 

\paragraph{Notations.} Throughout this paper, the following notations are used. Let $A$ be a generic matrix. Then $A^\top$ denotes its transpose. Further, $\mathbb{I}_{d\times d}$ denotes the $d\times d$-dimensional identity matrix, $\mathbb{O}_{d\times m}$ the $d \times m$-dimensional zero matrix, $0_d$ the $d$-dimensional vector of zeros, and $\textrm{diag}[a_1,\ldots,a_d]$ the $d\times d$-dimensional diagonal matrix with diagonal entries $a_1,\ldots,a_d$. Moreover, $\mathcal{N}_d(\mu,C)$ denotes the $d$-dimensional normal distribution with $d$-dimensional mean vector $\mu$ and $d\times d$-dimensional covariance matrix $C$. In addition, $|\cdot|$ denotes the absolute value or matrix norm.   Finally, parentheses~$(\cdot)$ and square brackets $[\cdot]$ denote vectors and matrices, respectively.


\section{Linear SDEs with additive noise}
\label{sec2:LinearSDEs}

Consider the $d$-dimensional linear SDE of additive noise type 
\begin{equation}\label{eq:LinearSDE}
    dX(t)=AX(t)dt+\Sigma dW(t), \qquad X(0)=X_0, 
\end{equation}
where $t\in [0,T]$, $T>0$, $A \in \mathbb{R}^{d\times d}$, $\Sigma\in \mathbb{R}^{d\times m}$, and $W(t)=(W^1(t),\ldots,W^m(t))^\top$, $t\in[0,T]$, is an $m$-dimensional standard Wiener process on the filtered probability space $(\Omega,\mathcal{F},(\mathcal{F}(t))_{t\in [0,T]},\mathbb{P})$, where $\mathcal{F}(t)=\sigma\{ W(s), \ 0<s<t \}$ is the sigma algebra generated by the Wiener process. The initial value $X_0$ is either deterministic, that is $X_0 \in \mathbb{R}^d$, or an $\mathbb{R}^d$-valued random variable which is independent of $(W(t))_{t\in [0,T]}$ and satisfies $\mathbb{E}\bigl[\norm{X_0}^{2}\bigr]<\infty$. 

It is well known (see, e.g., \cite{Mao2011}) that the unique strong solution to SDE \eqref{eq:LinearSDE} is, for all $t \in [0,T]$, given by
\begin{equation}\label{eq:solutionLinearSDE}
    X(t)=e^{At}X_0 + \eta(t), \qquad \eta(t) := \int\limits_0^t e^{A(t-s)} \Sigma \ dW(s),
\end{equation}
where 
\begin{equation}\label{eq:solutionLinearSDE_covMatrix}
    \eta(t) \sim \mathcal{N}_d\left( 0_d, C_\eta(t) \right), \qquad C_\eta(t):=\textrm{Cov}\left( \eta(t) \right) = \int\limits_0^t e^{A(t-s)} \Sigma\Sigma^\top \left( e^{A(t-s)} \right)^\top  ds \in \mathbb{R}^{d\times d}.
\end{equation}

\begin{remark}\label{rem:C_eta_full_rank}
    In the following, we assume that the covariance matrix $C_\eta(t)$ defined in \eqref{eq:solutionLinearSDE_covMatrix} is positive definite for every $t \in (0,T]$. This ensures that the $d$-dimensional normal distribution in~\eqref{eq:solutionLinearSDE_covMatrix} is non-degenerate, and therefore admits a density with respect to the Lebesgue measure. Depending on the noise structure of SDE \eqref{eq:LinearSDE}, two relevant  situations arise in this regard. In the \textbf{elliptic} case, the diffusion matrix $\Sigma\Sigma^\top$ is of full rank, which readily implies that $C_\eta(t)$ is of full rank for every $t \in (0,T]$. If $\Sigma\Sigma^\top$ is not of full rank, the drift matrix $A$ is such that  SDE~\eqref{eq:LinearSDE} is \textbf{hypoelliptic}. In this case, although the noise does not directly act on all components, it is propagated through the whole system via the drift coefficient, so that the covariance matrix $C_\eta(t)$ is nevertheless of full rank for every $t \in (0,T]$. 
\end{remark}


\paragraph{Standard exact simulation procedure.}

Let $0=t_0<t_1 < \ldots < t_n=T$ be an equidistant discretisation of the time interval $[0,T]$ with time step $\Delta=t_i-t_{i-1}$, for $i=1,\ldots,n$.

A realization of the process $(X(t))_{t\in [0,T]}$ in \eqref{eq:solutionLinearSDE} at the grid points can be obtained recursively via
\begin{equation*}
    X(t_i)=e^{A\Delta} X(t_{i-1}) + \eta_i(\Delta), \qquad i=1,\ldots,n,
\end{equation*}
where, similar as in \eqref{eq:solutionLinearSDE}-\eqref{eq:solutionLinearSDE_covMatrix},
\begin{equation}\label{eq:eta_i}
    \eta_i(\Delta):=\int\limits_{t_{i-1}}^{t_i} e^{A(t_i-s)}\Sigma \ dW(s) \sim \mathcal{N}_d(0_d,C_\eta(\Delta)), \qquad i=1,\ldots,n,
\end{equation}
with $C_\eta(\Delta)$ defined in \eqref{eq:solutionLinearSDE_covMatrix}. Here, the stochastic integrals $\eta_i(\Delta)$, $i=1,\ldots,n$, are independent. A~standard simulation procedure therefore consists of sampling realizations of $\eta_i(\Delta)$ from the $d$-dimensional normal distribution $\mathcal{N}_d(0_d,C_\eta(\Delta))$, using the unique Cholesky factorization of the covariance matrix $C_\eta(\Delta)$. 

Although this procedure generates an exact realization of the solution to the linear SDE \eqref{eq:LinearSDE}, in the sense that it preserves its distribution at the grid points, it does not explicitly relate the simulated path to the underlying realization of the driving Wiener process. 


\subsection{Wiener-path-based exact simulation procedure (according to \cite{Cruz2020})}
\label{sec2.1:LinearSDEs_WPB_exact}


\paragraph{Joint covariance structure.}

In \cite{Cruz2020}, the authors showed that paths of the process $(X(t))_{t \in [0,T]}$ in~\eqref{eq:solutionLinearSDE} can be related to their associated Wiener paths, by exploiting the joint covariance structure of the stochastic integral $\eta(t)$ defined in~\eqref{eq:solutionLinearSDE} and the Wiener increment 
\begin{equation}\label{eq:xi_t}
    \xi(t) 
    := 
    \int\limits_0^t \mathbb{I}_{m\times m} \ dW(s),
\end{equation}
where
\begin{equation}\label{eq:xi_t_covMatrix}
    \xi(t) \sim \mathcal{N}_m\left( 0_m, C_\xi(t) \right), \qquad C_\xi(t):=\textrm{Cov}\left( \xi(t) \right) = \int\limits_0^t \mathbb{I}_{m\times m} \ ds = t \mathbb{I}_{m\times m} \in \mathbb{R}^{m\times m}.
\end{equation}

Following \cite{Cruz2020}, we derive the $(m+d)\times(m+d)$-dimensional joint covariance matrix of $\xi(t)$  and $\eta(t)$, denoted by $C_{\xi;\eta}(t)$. Using It\^o's isometry, it is given by
\begin{align*}\label{eq:standard_expr_covMatrix}
    C_{\xi;\eta}(t) &:= \textrm{Cov}\Bigl( \begin{pmatrix}
        \xi(t) \\
        \eta(t)
    \end{pmatrix}  \Bigr) \notag  \\ 
    &= \textrm{Cov}\left( \int\limits_0^t  e^{  \begin{bmatrix}
        \mathbb{O}_{m \times m} & \mathbb{O}_{m\times d} \\
        \mathbb{O}_{d\times m} & A
    \end{bmatrix} (t-s) }  \begin{bmatrix}
        \mathbb{I}_{m\times m} \\
        \Sigma
    \end{bmatrix}  dW(s) \right) \notag  \\[0.1em]
    &= \int\limits_0^t  e^{ \begin{bmatrix}
        \mathbb{O}_{m \times m} & \mathbb{O}_{m\times d} \\
        \mathbb{O}_{d\times m} & A
    \end{bmatrix} (t-s) } \begin{bmatrix}
        \mathbb{I}_{m\times m} \\
        \Sigma
    \end{bmatrix} \begin{bmatrix}
        \mathbb{I}_{m\times m} & 
        \Sigma^\top
    \end{bmatrix}  e^{ \begin{bmatrix}
        \mathbb{O}_{m \times m} & \mathbb{O}_{m\times d} \\[0.1em]
        \mathbb{O}_{d\times m} & A
    \end{bmatrix}^\top (t-s) }  ds \\
    &= \int\limits_0^t  \begin{bmatrix}
        \mathbb{I}_{m\times m} & \Sigma^\top \left( e^{A(t-s)} \right)^\top \\
        e^{A(t-s)} \Sigma & e^{A(t-s)} \Sigma\Sigma^\top \left( e^{A(t-s)} \right)^\top
    \end{bmatrix}  ds. \notag  
\end{align*}
Therefore, it can  be expressed as
\begin{equation}\label{eq:C_xi_eta}
    C_{\xi;\eta}(t)=\begin{bmatrix}
        C_\xi(t) & C(t)^\top \\
        C(t) & C_\eta (t)
    \end{bmatrix},
\end{equation}
where $C_\eta(t)$ and $C_\xi(t)$  are as in \eqref{eq:solutionLinearSDE_covMatrix} and \eqref{eq:xi_t_covMatrix}, respectively, and the cross-covariance 
\begin{equation}\label{eq:covMatrix_C}
    C(t):= \textrm{Cov}(\eta(t),\xi(t)) = \int\limits_0^t e^{A(t-s)}\Sigma \ ds \ \in \mathbb{R}^{d\times m}.
\end{equation}

\begin{remark}
    In addition to the non-degeneracy of $C_\eta(t)$ (cf. Remark \ref{rem:C_eta_full_rank}), we assume that the joint covariance matrix $C_{\xi;\eta}(t)$ is positive definite for every $t \in \mathbb(0,T]$. This ensures that the joint  $(m+d)$-dimensional normal distribution of $(\xi(t),\eta(t))^\top$ is non-degenerate and that $C_{\xi;\eta}(t)$ admits a unique Cholesky factorization, which is used in the algorithm introduced below.
\end{remark}


\paragraph{Wiener-path-based exact simulation algorithm.}

Using the joint covariance matrix $C_{\xi;\eta}(\Delta)$ from~\eqref{eq:C_xi_eta}, one can sample joint realizations of the random vector 
\begin{equation*}
    \begin{pmatrix}
        \xi_i(\Delta) \\
        \eta_i(\Delta)
    \end{pmatrix},
\end{equation*}
where, similar as in \eqref{eq:xi_t}-\eqref{eq:xi_t_covMatrix},
\begin{equation}\label{eq:xi_i}
    \xi_i(\Delta):=\int\limits_{t_{i-1}}^{t_i} \mathbb{I}_{m\times m} \ dW(s) \sim \mathcal{N}_m(0_m,C_\xi(\Delta)), \qquad i=1,\ldots,n,
\end{equation}
are independent, and $C_\xi(\Delta)$ is defined in \eqref{eq:xi_t_covMatrix}. Specifically, these realizations are obtained by sampling from the joint $(m+d)$-dimensional normal distribution $\mathcal{N}_{m+d}(0_{m+d},C_{\xi;\eta}(\Delta))$, using the unique Cholesky factorization of $C_{\xi;\eta}(\Delta)$.

This procedure is summarized in Algorithm \ref{alg:Exact_WPB} and generates exact paths of $(X(t))_{t \in [0,T]}$ in \eqref{eq:solutionLinearSDE} that are coherent with their corresponding underlying Wiener paths.

\begin{algorithm}[t]
	\caption{\textbf{Wiener-path-based} exact simulation (\textbf{WPB-exact}, cf. \cite{Cruz2020}) 
		\ \\ \textbf{Input:} Initial value $X_0$, time horizon $T$, 
  step-size $\Delta$, exponential matrix~$e^{A \Delta}$ and joint covariance matrix $C_{\xi;\eta}(\Delta)$. 
  \ \\
		\textbf{Output:} Path of the linear SDE \eqref{eq:LinearSDE} and corresponding Wiener path. 
        }
	\label{alg:Exact_WPB}
    \begin{spacing}{1.1}
	\begin{algorithmic}[1]
        \State Set $n=T/\Delta$.
		\State Set ${X}(t_0)=X_0$ and $W(t_0)=0_m$. 
        \State Determine $L$ such that $L L^\top=C_{\xi;\eta}(\Delta)$.
		\For{$i = 1:n$}
        \State Generate $\mathbf{z}=(z_1,\ldots,z_{m+d})^\top$ from $\mathcal{N}_{m+d}\bigl( 0_{m+d}, \mathbb{I}_{(m+d)\times (m+d)} \bigr)$.
        \State Set $(\mathbf{y},\mathbf{x})^\top=(y_1,\ldots,y_m,x_1,\ldots,x_d)^\top=L\mathbf{z}$.
		\State Set $X(t_i)=e^{A\Delta}X(t_{i-1})+\mathbf{x}$.  
        \State Set $W(t_i)=W(t_{i-1})+\mathbf{y}$.
		\EndFor
		\State Return ${X}(t_i)$, 
        $i=0,\ldots,n$, and $W(t_i)$, $i=0,\ldots,n$. 
	\end{algorithmic}
    \end{spacing}
\end{algorithm}

\begin{remark}
    Algorithm \ref{alg:Exact_WPB} requires sampling from an $(m+d)$-dimensional standard normal distribution. In the special case of $d=m$, it is natural to ask whether a more efficient procedure can be obtained by sampling only a realization $\mathbf{z}=(z_1,\ldots,z_d)^\top$ of a $d$-dimensional standard normal distribution and using this same random source to propagate both the solution path and the Wiener path. In Section \ref{appA}, we propose such an alternative Wiener-path-based simulation procedure, which reduces the number of required random variables while providing a controlled approximation of the coupling between the solution path and its driving Wiener path.  
\end{remark}


\subsection{Wiener-path-based exact simulation on different time grids}
\label{sec2.2:LinearSDEs_WPB_exact_grids}

In this section, we extend the procedure of \cite{Cruz2020} (cf. Algorithm \ref{alg:Exact_WPB}) to generate a realization of SDE~\eqref{eq:LinearSDE} simultaneously on two different time grids. 

Consider two different time step-sizes $\Delta_1<\Delta_2$, where $\Delta_2=l\Delta_1$ for some $l\in \mathbb{N}$. Then, we aim to generate a path~${X}_{\Delta_1}(t_i)$, $t_i=i\Delta_1$, $i=0,\ldots,n_1$, where $t_{n_1}=T$, together with the corresponding path ${X}_{\Delta_2}(t_j)$ on the coarser time grid $t_j=j\Delta_2$, $j=0,\ldots,n_2$, where $t_{n_2}=T$, such that both approximations are based on the same underlying Wiener path. More precisely, we require to obtain
\begin{equation*}
    X_{\Delta_2}(t_j)=X_{\Delta_1}(t_{lj}), \quad W_{\Delta_2}(t_j)=W_{\Delta_1}(t_{lj}), \quad j=1,\ldots,n_2.
\end{equation*}

\paragraph{Preliminary analysis.}

To achieve this, we first consider one step of size $\Delta_2$ from a time point $t$ to $t+\Delta_2$ of the solution and Wiener process on the coarse grid. Similar as in \eqref{eq:solutionLinearSDE}-\eqref{eq:solutionLinearSDE_covMatrix} and \eqref{eq:xi_t}-\eqref{eq:xi_t_covMatrix}, they are given by
\begin{equation}\label{eq:eta_Delta2}
    X(t+\Delta_2)=e^{A\Delta_2}X(t) + \eta(\Delta_2), 
\end{equation}
where
\begin{equation*}
    \eta(\Delta_2)=\int\limits_t^{t+\Delta_2} e^{A(t+\Delta_2-s)} \Sigma \ dW(s) \sim \mathcal{N}_d(0_d,C_\eta(\Delta_2)),
\end{equation*}
and by
\begin{equation}\label{eq:xi_Delta2}
    W(t+\Delta_2)=W(t) + \xi(\Delta_2),
\end{equation}
where
\begin{equation*}
    \xi(\Delta_2)=\int\limits_{t}^{t+\Delta_2} \mathbb{I}_{m\times m} \ dW(s) 
    \sim \mathcal{N}_m(0_m,C_\xi(\Delta_2)).
\end{equation*}
Moreover, the joint covariance matrix of $(\xi(\Delta_2),\eta(\Delta_2))^\top$ is given by $C_{\xi;\eta}(\Delta_2)$ as in \eqref{eq:C_xi_eta}.

The key question is how to construct the random vector $(\xi(\Delta_2),\eta(\Delta_2))^\top$
based on the $l$ corresponding random vectors, denoted by $(\xi_k(\Delta_1),\eta_k(\Delta_1))^\top$, $k=1,\ldots,l$, 
used to propagate the solution and Wiener process on the fine grid over the same time increment?

To this end, note that similar as above, $l$ steps on the fine grid, from $t$ to $t+l\Delta_1$, can be expressed via
\begin{equation}\label{eq:eta_k_Delta1}
    X(t+l\Delta_1) = e^{Al\Delta_1}X(t)+\sum\limits_{k=1}^{l} e^{A(l-k)\Delta_1} \eta_k(\Delta_1),
\end{equation}
where
\begin{equation*}
    \eta_k(\Delta_1)=\int\limits_{t+(k-1)\Delta_1}^{t+k\Delta_1} e^{A(t+k\Delta_1-s)}\Sigma \ dW(s) \sim \mathcal{N}_d(0_d,C_\eta(\Delta_1)), \quad k=1,\ldots,l.
\end{equation*}
Moreover, for the Wiener process we have
\begin{equation}\label{eq:xi_k_Delta1}
    W(t+l\Delta_1)=W(t)+\sum\limits_{k=1}^l \xi_k(\Delta_1),
\end{equation}
where
\begin{equation*}
    \xi_k(\Delta_1)=\int\limits_{t+(k-1)\Delta_1}^{t+k\Delta_1} \mathbb{I}_{m\times m} \ dW(s) \sim \mathcal{N}_m(0_m,C_\xi(\Delta_1)), \quad k=1,\ldots,l.
\end{equation*}
Now, recalling that $\Delta_2=l\Delta_1$ and comparing \eqref{eq:eta_Delta2} and \eqref{eq:eta_k_Delta1} as well as \eqref{eq:xi_Delta2} and \eqref{eq:xi_k_Delta1}, suggests defining
\begin{equation}\label{eq:eta_xi_coarse}
    \eta_{\Delta_2}^{\textrm{coarse}}:=\sum\limits_{k=1}^l e^{A(l-k)\Delta_1} \eta_k(\Delta_1), \quad \xi_{\Delta_2}^{\textrm{coarse}}:=\sum\limits_{k=1}^l \xi_k(\Delta_1).
\end{equation}

In the following proposition, we show that, indeed,
\begin{equation*}
    \eta_{\Delta_2}^{\textrm{coarse}}=\eta(\Delta_2), \qquad \xi_{\Delta_2}^{\textrm{coarse}}=\xi(\Delta_2).
\end{equation*} 

\begin{proposition}\label{prop:coarse_eta_xi}
    Let $\Delta_1<\Delta_2$ such that $\Delta_2=l\Delta_1$, for some $l \in \mathbb{N}$. Moreover, let $\eta_{\Delta_2}^{\textrm{coarse}}$ and $\xi_{\Delta_2}^{\textrm{coarse}}$ be as in \eqref{eq:eta_xi_coarse}. Then, they equal $\eta(\Delta_2)$ and $\xi(\Delta_2)$ as in \eqref{eq:eta_Delta2} and \eqref{eq:xi_Delta2}, respectively.
\end{proposition}

\begin{proof}
    It holds that
    \begin{samepage}
        \begin{align*}
        \eta_{\Delta_2}^{\textrm{coarse}} &= \sum\limits_{k=1}^l e^{A(l-k)\Delta_1} \eta_k(\Delta_1) \\
        &= \sum\limits_{k=1}^l e^{A(l-k)\Delta_1} \int\limits_{t+(k-1)\Delta_1}^{t+k\Delta_1} e^{A(t+k\Delta_1-s)}\Sigma \ dW(s) \\
        &= \sum\limits_{k=1}^l \ \int\limits_{t+(k-1)\Delta_1}^{t+k\Delta_1} e^{A(t+l\Delta_1-s)} \Sigma \ dW(s) \\
        &=\int\limits_t^{t+l\Delta_1} e^{A(t+l\Delta_1-s)} \Sigma \ dW(s) \\
        &=\int\limits_t^{t+\Delta_2} e^{A(t+\Delta_2-s)} \Sigma \ dW(s) 
        =\eta(\Delta_2).
        \end{align*}
    \end{samepage}

    Moreover, one has
    \begin{align*}
        \xi_{\Delta_2}^{\textrm{coarse}} &= \sum\limits_{k=1}^l \xi_k(\Delta_1) \\
        &= \sum\limits_{k=1}^l \ \int\limits_{t+(k-1)\Delta_1}^{t+k\Delta_1} \mathbb{I}_{m\times m} \ dW(s) \\
        &= \int\limits_{t}^{t+l\Delta_1} \mathbb{I}_{m\times m} \ dW(s) \\
        &= \int\limits_{t}^{t+\Delta_2} \mathbb{I}_{m\times m} \ dW(s) = \xi(\Delta_2).
    \end{align*}

\end{proof}

\noindent Proposition \ref{prop:coarse_eta_xi} shows that the coarse-grid random variables  $(\xi(\Delta_2),\eta(\Delta_2))^\top$ can indeed be constructed from the fine-grid random vectors $(\xi_k(\Delta_1),\eta_k(\Delta))^\top$, $k=1,\ldots,l$, via \eqref{eq:eta_xi_coarse}. In particular, this construction preserves both the marginal distributions of $\xi(\Delta_2)$ and $\eta(\Delta_2)$ and their joint covariance structure $C_{\xi;\eta}(\Delta_2)$.

\paragraph{Wiener-path-based exact simulation algorithm on different time grids.}

Using the previous results, the corresponding algorithm is obtained as follows. 
For each coarse-grid step $j=1,\ldots,n_2$, the fine-grid path $X_{\Delta_1}$ and its associated Wiener path $W_{\Delta_1}$ are first propagated over the $l$ corresponding fine-grid steps as in Algorithm \ref{alg:Exact_WPB}. This is done by sampling independent realizations $\mathbf{z}_k$, $k=1,\ldots,l$, from the $(m+d)$-dimensional standard normal distribution, and setting $(\mathbf{y}_k,\mathbf{x}_k)^\top=L_{\Delta_1}\mathbf{z}_k$, where $L_{\Delta_1} L_{\Delta_1}^\top=C_{\xi;\eta}(\Delta_1)$. The coarse-grid increments are then obtained by
\begin{equation*}
    \mathbf{x}=\sum\limits_{k=1}^l e^{A(l-k)\Delta_1} \mathbf{x}_k, \quad \mathbf{y}=\sum\limits_{k=1}^l \mathbf{y}_k,
\end{equation*}
and are subsequently used to propagate the coarse-grid path $X_{\Delta_2}$ and the corresponding Wiener path $W_{\Delta_2}$. The resulting procedure is summarized in Algorithm \ref{alg:Exact_WPB_twoDeltas}.

\begin{algorithm}[t]
	\caption{\textbf{Wiener-path-based} exact simulation on different time grids (\textbf{WPB-exact-grids}) 
		\ \\ \textbf{Input:} Initial value $X_0$, time horizon $T$, 
  step-sizes $\Delta_1$ and $\Delta_2$, drift matrix $A$,  exponential matrices $e^{A \Delta_1}$ and $e^{A \Delta_2}$ and joint covariance matrix $C_{\xi;\eta}(\Delta_1)$.
  \ \\
		\textbf{Output:} Path of the linear SDE \eqref{eq:LinearSDE} and corresponding Wiener path, generated on two different time grids. 
        }
	\label{alg:Exact_WPB_twoDeltas}
    \begin{spacing}{1.1}
	\begin{algorithmic}[1]
        \State Set $n_2=T/\Delta_2$, $l=\Delta_2/\Delta_1$ and $i=0$.
		\State Set ${X}_{\Delta_1}(t_0)={X}_{\Delta_2}(t_0)=X_0$ and $W_{\Delta_1}(t_0)=W_{\Delta_2}(t_0)=0_m$. 
        \State Determine $L_{\Delta_1}$ such that $L_{\Delta_1} L_{\Delta_1}^\top=C_{\xi;\eta}(\Delta_1)$.
		\For{$j = 1:n_2$}
        \For{$k = 1:l$}
        \State Generate $\mathbf{z}_k$ from $\mathcal{N}_{m+d}\bigl( 0_{m+d}, \mathbb{I}_{(m+d)\times (m+d)} \bigr)$ and set $(\mathbf{y}_k,\mathbf{x}_k)^\top=L_{\Delta_1}\mathbf{z}_k$.
        \State Set $i=i+1$ (or $i=(j-1)l+k$).
		\State Set $X_{\Delta_1}(t_i)=e^{A\Delta_1}X_{\Delta_1}(t_{i-1})+\mathbf{x}_k$.  
        \State Set $W_{\Delta_1}(t_i)=W_{\Delta_1}(t_{i-1})+\mathbf{y}_k$.
		\EndFor
        \State Set $\mathbf{y}=\sum\limits_{k=1}^{l}\mathbf{y}_k$.
        \State Set $\mathbf{x}=\sum\limits_{k=1}^{l} e^{A(l-k)\Delta_1} \mathbf{x}_k$.
		\State Set $X_{\Delta_2}(t_j)=e^{A\Delta_2}X_{\Delta_2}(t_{j-1})+\mathbf{x}$.  
        \State Set $W_{\Delta_2}(t_j)=W_{\Delta_2}(t_{j-1})+\mathbf{y}$.
        \EndFor
		\State Return ${X}_{\Delta_1}(t_i)$ and $W_{\Delta_1}(t_i)$, $i=0,\ldots,n_1$, and ${X}_{\Delta_2}(t_j)$ and $W_{\Delta_2}(t_j)$, $j=0,\ldots,n_2$. 
	\end{algorithmic}
    \end{spacing}
\end{algorithm}


\section{Illustration of the algorithms on the simple stochastic harmonic oscillator}
\label{sec3:simple_SHO}

In this section, the \textit{WPB-exact} algorithm by \cite{Cruz2020} (Algorithm \ref{alg:Exact_WPB}) and the proposed \textit{WPB-exact-grids} algorithm (Algorithm \ref{alg:Exact_WPB_twoDeltas}) are illustrated on the simple stochastic harmonic oscillator (SHO), a well-known stochastic model from classical mechanics.

\paragraph{Simple stochastic harmonic oscillator.}

The simple SHO (see, e.g., \cite{Cohen2012,SENOSIAIN2019,Strommen2004}) is described by the following ($d=2$)-dimensional linear SDE with additive noise:
\begin{equation}\label{eq:SimpleSHO_elliptic}
	d \underbrace{\begin{pmatrix}
			Q(t) \\
			P(t) 
	\end{pmatrix}}_{=X(t)}
	=
	\underbrace{\begin{pmatrix}
			P(t) \\
			-Q(t)
	\end{pmatrix}}_{=AX(t)} dt \ + \
	\underbrace{\begin{bmatrix}
		\sigma_1 &	\sigma_3 \\ 
        \sigma_4 & \sigma_2
	\end{bmatrix}}_{=\Sigma} dW(t), \qquad X(0)=X_0, 
\end{equation}
where $t\in[0,T]$, $W(t)=(W^1(t),W^2(t))^\top$, $t\in [0,T]$, is an $(m=2)$-dimensional standard Wiener process, $\sigma_k$, $k=1,\ldots,4$, denote the noise intensity parameters, and the drift matrix
\begin{equation*}
    A=\begin{bmatrix}
        0 & 1 \\
        -1 & 0
    \end{bmatrix}.
\end{equation*}

\begin{remark}
    Depending on the structure of the diffusion matrix $\Sigma$, SDE \eqref{eq:SimpleSHO_elliptic} may be either elliptic or hypoelliptic. If $\Sigma$ has full rank, the system is elliptic. A standard example is obtained by choosing a diagonal diffusion matrix, i.e., $\sigma_3=\sigma_4=0$ and $\sigma_1,\sigma_2>0$. In contrast, if the noise acts only on the second component, that is,  $\sigma_1=\sigma_3=\sigma_4=0$ and $\sigma_2>0$, the resulting SDE is hypoelliptic. In this case, the driving noise can equivalently be represented by an $(m=1)$-dimensional Wiener process, and the joint covariance matrix  $C_{\xi;\eta}(t)$ introduced below reduces to dimension $3\times 3$. Algorithms \ref{alg:Exact_WPB} and \ref{alg:Exact_WPB_twoDeltas} apply analogously to this case, yielding Wiener-path-based realizations with the same properties as in the elliptic setting. Throughout this section, we consider the classical elliptic case with diagonal diffusion matrix. Additional choices of the diffusion matrix $\Sigma$ are discussed in Section \ref{appA}.   
\end{remark}

\paragraph{Matrix exponential and joint covariance matrix.}

In the following, we derive closed-form expressions for the matrix exponential $e^{A\Delta}$ and joint covariance matrix $C_{\xi;\eta}(\Delta)$, which serve as inputs for Algorithm \ref{alg:Exact_WPB} and Algorithm \ref{alg:Exact_WPB_twoDeltas}.

For $t>0$, the matrix exponential is given by
\begin{equation}\label{eq:ExpMat_SHO}
    e^{At}=\begin{bmatrix}
        \cos(t) & \sin(t) \\
        -\sin(t) & \cos(t)
    \end{bmatrix}.
\end{equation}
Moreover, for $t>0$, the joint covariance matrix takes the form
\begin{equation}\label{eq:JointCovMat_SHO}
    C_{\xi;\eta}(t)=\begin{bmatrix}
        C_\xi(t) & C(t)^\top \\
        C(t) & C_\eta(t)
    \end{bmatrix},
\end{equation}
where
\begin{align*}
    C_\xi(t)&=t\mathbb{I}_{2\times 2}, \\[1em]
    C_\eta(t)
    &=\frac{1}{2}\begin{bmatrix}
        c_{\eta;11}(t) & c_{\eta;12(t)} \\
        c_{\eta;12(t)} & c_{\eta;22(t)}
    \end{bmatrix}, \quad \text{with} \\[1em]
    & c_{\eta;11}(t) = \left( \sigma_1^2+\sigma_2^2+\sigma_3^2+\sigma_4^2 \right)t-\left( \sigma_2^2+\sigma_4^2-\sigma_1^2-\sigma_3^2 \right)\cos(t)\sin(t)+2\left( \sigma_2\sigma_3+\sigma_1\sigma_4 \right)\sin^2(t), \\[1em]
    & c_{\eta;12}(t) = \left( \sigma_2^2+\sigma_4^2-\sigma_1^2-\sigma_3^2 \right)\sin^2(t)+\left( \sigma_2\sigma_3+\sigma_1\sigma_4 \right)\sin(2t), \\[1em]
    & c_{\eta;22}(t) = \left( \sigma_1^2+\sigma_2^2+\sigma_3^2+\sigma_4^2 \right)t+\left( \sigma_2^2+\sigma_4^2-\sigma_1^2-\sigma_3^2 \right)\cos(t)\sin(t)-2\left( \sigma_2\sigma_3+\sigma_1\sigma_4 \right)\sin^2(t), \\[1em]
    C(t) &= \begin{bmatrix}
    \sigma_1 \sin(t)+\sigma_4\left( 1-\cos(t) \right) & \sigma_2\left( 1-\cos(t) \right)+\sigma_3\sin(t) \\[0.5em]
    \sigma_1\left( \cos(t)-1 \right) + \sigma_4\sin(t) & \sigma_2 \sin(t) + \sigma_3 \left( \cos(t)-1 \right)
    \end{bmatrix}. 
\end{align*}
Note that, while $C_\xi(t)$ and $C_\eta(t)$ are covariance matrices and hence symmetric, the cross-covariance matrix $C(t)$ is generally not symmetric.

\paragraph{Numerical experiments.}

Algorithm \ref{alg:Exact_WPB} is now applied to the simple SHO \eqref{eq:SimpleSHO_elliptic} using the input quantities $e^{A\Delta}$ and $C_{\xi;\eta}(\Delta)$ from \eqref{eq:ExpMat_SHO} and \eqref{eq:JointCovMat_SHO}, respectively. The initial value is set to $X_0=(0,0)^\top$, the time horizon to $T=50$, and the step-size to $\Delta=10^{-4}$. The noise parameters are chosen as $\sigma_1=1$, $\sigma_2=2$, and $\sigma_3=\sigma_4=0$.

\begin{figure}[t]
	\begin{centering}
		\includegraphics[width=1.0\textwidth]{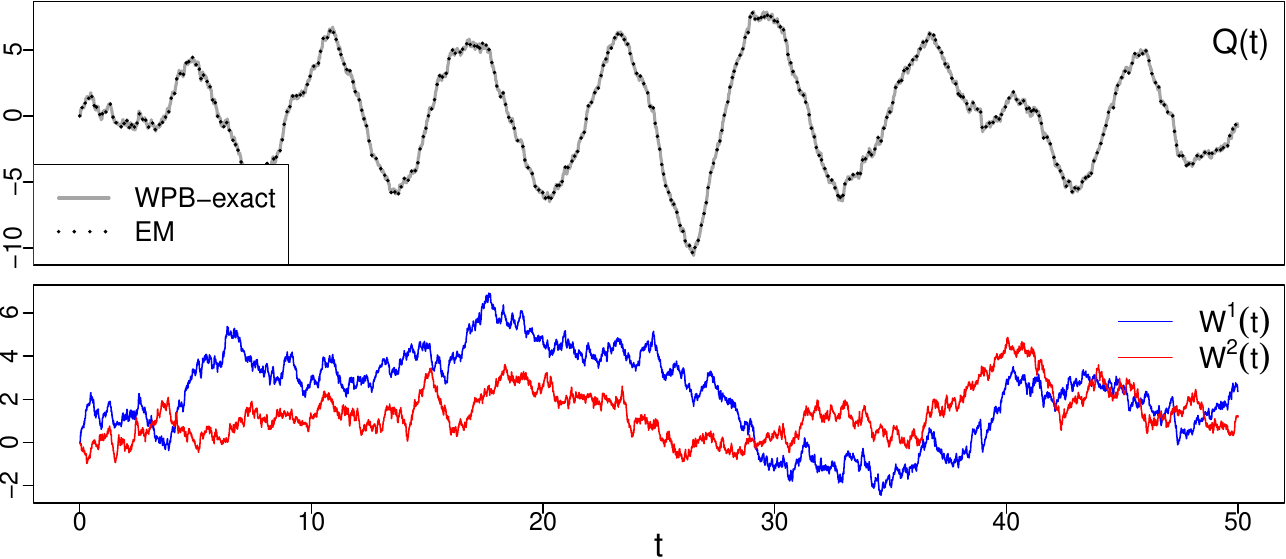}	
		\caption{Path of the $Q$-component of the simple SHO \eqref{eq:SimpleSHO_elliptic}. Top panel: realization $Q(t_i)$ (gray solid line) obtained via Algorithm \ref{alg:Exact_WPB}, compared with the Euler--Maruyama realization $\widetilde{Q}^{\textrm{EM}}(t_i)$ (black dotted line) obtained using the same underlying Wiener path. Bottom panel: corresponding Wiener path $W(t_i)=(W^1(t_i),W^2(t_i))^\top$, $i=0,\ldots,n$.}
		\label{fig:simpleSHO_elliptic}
	\end{centering}
\end{figure}

Figure \ref{fig:simpleSHO_elliptic} shows the resulting realization of the $Q$-component, $Q(t_i)$, $i=0,\ldots,n$, generated by Algorithm \ref{alg:Exact_WPB} (top panel; gray solid line), together with the corresponding underlying Wiener path $W(t_i)=(W^1(t_i),W^2(t_i))^\top$, $i=0,\ldots,n$ (bottom panel; blue and red solid lines).

Furthermore, the Wiener path generated by Algorithm \ref{alg:Exact_WPB} is used to construct a corresponding approximation via the standard Euler--Maruyama (EM) scheme \cite{Kloeden1992}. The resulting approximate path $\widetilde{Q}^{\textrm{EM}}(t_i)$, $i=0,\ldots,n$, is also shown in the top panel of Figure \ref{fig:simpleSHO_elliptic} (black dotted line). As expected, due to the small step size $\Delta=10^{-4}$, the EM approximation closely matches the path $Q(t_i)$, $i=0,\ldots,n$, generated by Algorithm \ref{alg:Exact_WPB}, the two realizations being  visually indistinguishable on the scale of the figure.

We now apply Algorithm \ref{alg:Exact_WPB_twoDeltas} to the simple SHO \eqref{eq:SimpleSHO_elliptic}. The simulation setup is identical to the previous experiment, except for the step sizes, which are now set to $\Delta_1=10^{-4}$ and $\Delta_2=10^{-2}$. 

\begin{figure}[t]
	\begin{centering}
		\includegraphics[width=1.0\textwidth]{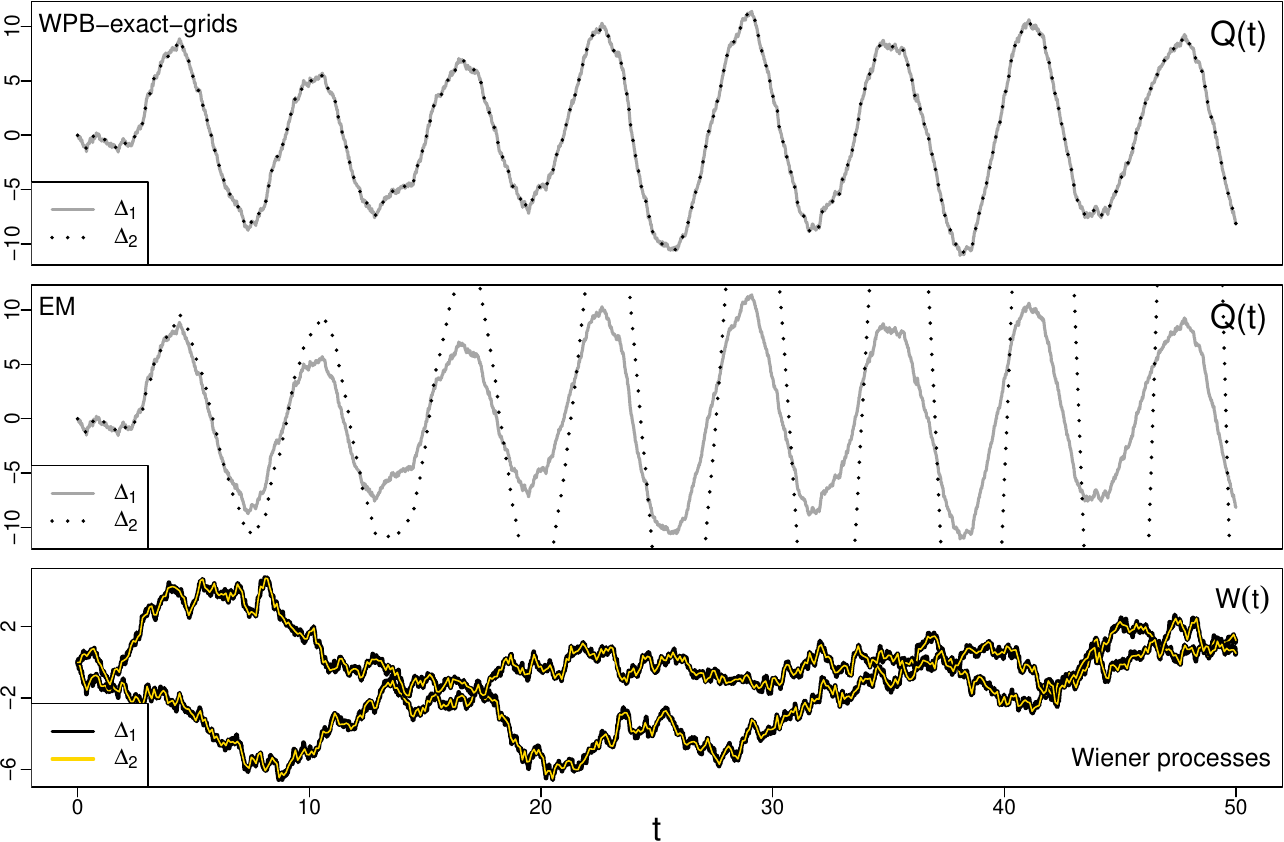}	
		\caption{Path of the $Q$-component of the simple SHO \eqref{eq:SimpleSHO_elliptic}. Top panel: fine- and coarse-grid paths  $Q_{\Delta_1}(t_i)$ (gray solid line) and $Q_{\Delta_2}(t_j)$ (black dotted line) obtained via Algorithm~\ref{alg:Exact_WPB_twoDeltas}. Middle panel: corresponding Euler--Maruyama realizations $\widetilde{Q}^{\textrm{EM}}_{\Delta_1}(t_i)$ (gray solid line) and $\widetilde{Q}^{\textrm{EM}}_{\Delta_2}(t_j)$ (black dotted line). Bottom panel: underlying Wiener paths $W_{\Delta_1}(t_i)$ (black lines) and $W_{\Delta_2}(t_j)$ (yellow lines).}
		\label{fig:simpleSHO_elliptic_twoDeltas}
	\end{centering}
\end{figure}

Figure \ref{fig:simpleSHO_elliptic_twoDeltas} shows the resulting paths 
\begin{equation*}
    Q_{\Delta_1}(t_i), \ i=0,\ldots,n_1, \quad \text{and} \quad Q_{\Delta_2}(t_j), \ j=0,\ldots,n_2, 
\end{equation*}
on the fine and coarse grids, respectively (top panel: gray solid and black dotted lines). The associated Wiener paths
\begin{equation*}
    W_{\Delta_1}(t_i)=(W^1_{\Delta_1}(t_i),W^2_{\Delta_1}(t_i))^\top, \quad i=0,\ldots,n_1,
\end{equation*}
and
\begin{equation*}
    W_{\Delta_2}(t_j)=(W^1_{\Delta_2}(t_j),W^2_{\Delta_2}(t_j))^\top, \quad j=0,\ldots,n_2,
\end{equation*}
are displayed in the bottom panel (black and yellow lines, respectively). 

It can be observed that the coarse-grid path $Q_{\Delta_2}(t_j)$ obtained with $\Delta_2=10^{-2}$ coincides with the fine-grid realization $Q_{\Delta_1}(t_i)$ obtained with $\Delta_1=10^{-4}$, and that also $W_{\Delta_2}(t_j)$ coincides with $W_{\Delta_1}(t_i)$. This agreement illustrates that Algorithm~\ref{alg:Exact_WPB_twoDeltas} generates the same discrete realization on both time grids, as expected from Proposition~\ref{prop:coarse_eta_xi}. 

In addition, the Wiener paths $W_{\Delta_1}(t_i)$ and $W_{\Delta_2}(t_j)$, returned by Algorithm \ref{alg:Exact_WPB_twoDeltas}, are used to generate the corresponding Euler--Maruyama approximations $\widetilde{Q}^{\textrm{EM}}_{\Delta_1}(t_i)$ and $\widetilde{Q}^{\textrm{EM}}_{\Delta_2}(t_j)$.  
These paths are shown in the middle panel of Figure \ref{fig:simpleSHO_elliptic_twoDeltas} (gray solid and black dotted line, respectively).

It can be observed that the fine-grid EM realization $\widetilde{Q}^{\textrm{EM}}_{\Delta_1}(t_i)$ closely approximates the exact path. However, the coarse-grid approximation $\widetilde{Q}^{\textrm{EM}}_{\Delta_2}(t_j)$ deviates significantly and fails to preserve the characteristic oscillatory dynamics of the simple SHO. This behaviour highlights a known limitation of the EM scheme when applied to oscillatory systems: although the method converges as the step size tends to zero, it may fail to reproduce the qualitative properties of the underlying dynamics for larger step sizes (see, e.g., \cite{Cohen2012,Strommen2004}).

\section{Semi-linear SDEs with additive noise}
\label{sec4:SemiLinearSDEs}

Consider the $d$-dimensional non-linear SDE of additive noise type
\begin{equation}\label{eq:semi-linearSDE}
    dX(t)=F(X(t))dt+\Sigma dW(t), \qquad X(0)=X_0,
\end{equation}
with $t \in [0,T]$, $T>0$, $\Sigma \in \mathbb{R}^{d\times m}$, $(W(t))_{t\in [0,T]}$ as for SDE \eqref{eq:LinearSDE}, and semi-linear drift coefficient given by
\begin{equation*}
    F(X(t)):=AX(t)+N(X(t)),
\end{equation*}
where $A\in\mathbb{R}^{d\times d}$ and $N\colon \mathbb{R}^d \to \mathbb{R}^d$.

In this section, we first recall numerical splitting methods for SDE \eqref{eq:semi-linearSDE} (Section \ref{sec4.1:Splitting_construction}). In Section~\ref{sec4.2:SemiLinearSDEs_WPB_splitting}, we then show how the joint covariance matrix $C_{\xi;\eta}(\Delta)$ from \eqref{eq:C_xi_eta} can be incorporated into these methods to generate approximate paths of \eqref{eq:semi-linearSDE} that are coherent with their underlying Wiener paths, thereby extending Algorithms \ref{alg:Exact_WPB} and \ref{alg:Exact_WPB_twoDeltas}. 

\begin{remark}\label{rem:N_assumptions}
    In agreement with the literature, we assume that the function $N$ is either globally Lipschitz continuous or satisfies a global one-sided Lipschitz condition and grows at most polynomially. In the first case, we assume that the initial condition $X_0$, independent of $(W(t))_{t\in [0,T]}$, satisfies $\mathbb{E}[\norm{X_0}^2]<\infty$, whereas in the latter case we require $\mathbb{E}[\norm{X_0}^{2p}]<\infty$ for all $p\geq 1$. Under these assumptions, SDE \eqref{eq:semi-linearSDE} admits a unique strong solution, which is defined for all $t\in[0,T]$ and has bounded moments 
    (see, e.g.,~\mbox{\cite{Alyushina1988,Hutzenthaler2012_2,Kelly2017,Mao2011,Tretyakov2012}}).

    Furthermore, as for linear SDEs (cf. Remark \ref{rem:C_eta_full_rank}), we assume that whenever $\Sigma\Sigma^\top$ is not of full rank, SDE \eqref{eq:semi-linearSDE} is hypoelliptic. This ensures that the noise is propagated through the whole system via the drift term, the solution  admitting a smooth density.
\end{remark}


\subsection{Construction of splitting methods}
\label{sec4.1:Splitting_construction}

Numerical splitting methods are typically constructed according to the following three steps: 
\begin{itemize}
	\item[(i)] Split the original (generally non-solvable) equation into exactly solvable subequations; 
	\item[(ii)] Derive the exact solution of each subequation; 
	\item[(iii)] Compose these exact solutions appropriately at each time step.
\end{itemize}

In the following, we recall known splitting methods for SDE \eqref{eq:semi-linearSDE} (see, e.g, \cite{Ableidinger2017,Buckwar2022,Pilipovic2024}).

\paragraph*{Step (i): Definition of subequations.}
We split SDE \eqref{eq:semi-linearSDE} into the following two subequations
\begin{align}
	d X^{[1]}(t)&=AX^{[1]}(t) dt + \Sigma dW(t), \qquad X^{[1]}(0)=X^{[1]}_0, 
	\label{SDE} \\
	dX^{[2]}(t)&=N(X^{[2]}(t))dt, \qquad X^{[2]}(0)=X^{[2]}_0. 
	\label{ODE}
\end{align} 

\paragraph*{Step (ii): Solutions of subequations.}
Subequation \eqref{SDE} is a linear SDE with additive noise, i.e. it is of type \eqref{eq:LinearSDE}. Therefore, its solution is given for all $t \in [0,T]$ by
\begin{equation}\label{eq:Exact_SDE}
  \varphi_t^{[1]}\bigl(X_0^{[1]}\bigr)=X^{[1]}(t)=e^{At}X_0^{[1]} + \eta(t), 
\end{equation}
where $\eta(t)$ is as in \eqref{eq:solutionLinearSDE} and has a $d$-dimensional normal distribution with covariance matrix $C_\eta(t)$ as in \eqref{eq:solutionLinearSDE_covMatrix}.

Subequation \eqref{ODE} is an ODE, which treats the non-linear term $N$. Due to the assumptions for~$N$ (see Remark \ref{rem:N_assumptions}), its unique solution exists and is defined on the entire time interval $[0,T]$. Denoting this solution by $f \colon \mathbb{R}^d \to \mathbb{R}^d$, we obtain
\begin{equation}\label{eq:Exact_ODE}
	\varphi_{t}^{[2]}(X^{[2]}_0)=X^{[2]}(t)=f(X^{[2]}_0;t). 
\end{equation}

\begin{remark}\label{rem:A_assumptions}
	Following \cite{Buckwar2022}, whenever $\Sigma\Sigma^\top$ is not of full rank, we choose the matrix $A$ in subequation~\eqref{SDE} such that it is hypoelliptic. As in Section \ref{sec2:LinearSDEs}, this guarantees that $C_\eta(t)$ is positive definite, and, consequently, that the corresponding $d$-dimensional normal distribution is non-degenerate. 
\end{remark}

\paragraph*{Step (iii): Composition of solutions.}

Similar to Section \ref{sec2:LinearSDEs}, let $0=t_0<t_1<\ldots < t_n=T$ be an equidistant  discretisation of the time interval $[0,T]$ with time step $\Delta=t_i-t_{i-1} \in (0,\Delta_{0}]$, $\Delta_{0} \in (0,1)$, for $i=1,\ldots,n$.

To obtain numerical approximations of the solution to SDE \eqref{eq:semi-linearSDE} at the grid points, the flows \eqref{eq:Exact_SDE} and \eqref{eq:Exact_ODE} of the two subequations have to be composed appropriately at each time step. Two common composition approaches are the Lie-Trotter (LT) \cite{Trotter1959} and the Strang (S) \cite{Strang1968} splitting. These yield the numerical schemes
\begin{equation}\label{eq:splitting}
    \begin{aligned}
	\widetilde{X}^{\textrm{LT}}(t_i)&=\left( \varphi_\Delta^{[1]} \circ \varphi_\Delta^{[2]} \right) \bigl(\widetilde{X}^{\textrm{LT}}(t_{i-1}) \bigr)
    =e^{A\Delta}f(\widetilde{X}^{\textrm{LT}}(t_{i-1});\Delta)+\eta_{i}(\Delta), \\ 
    \widetilde{X}^{\textrm{S}}(t_i)&=\left( \varphi_{\Delta/2}^{[2]} \circ \varphi_{\Delta}^{[1]} \circ \varphi_{\Delta/2}^{[2]} \right)\bigl(\widetilde{X}^{\textrm{S}}(t_{i-1}) \bigr)=f\left(e^{A\Delta}f\bigl(\widetilde{X}^{\textrm{S}}(t_{i-1});\Delta/2\bigr) + \eta_{i}(\Delta);\Delta/2 \right),
\end{aligned}
\end{equation}
for $i=1,\ldots,n$, where $\widetilde{X}^{\textrm{LT}}(t_0)=\widetilde{X}^{\textrm{S}}(t_0)=X_0$ and the $\eta_i(\Delta)$, $i=1,\ldots,n$, are as in \eqref{eq:eta_i}. 

\begin{remark}
    The numerical splitting schemes \eqref{eq:splitting} have been proved to be strongly convergent (in the mean-square sense) of order $1$ (see, e.g., \cite{Ableidinger2017,Buckwar2022,Pilipovic2024}).
\end{remark}


\subsection{Wiener-path-based splitting simulation algorithms}
\label{sec4.2:SemiLinearSDEs_WPB_splitting}

In this subsection, we discuss how the joint covariance matrix $C_{\xi;\eta}(\Delta)$ from \eqref{eq:C_xi_eta} (cf. Section \ref{sec2:LinearSDEs}) can be combined with the splitting methods \eqref{eq:splitting} to generate approximate paths of SDE \eqref{eq:semi-linearSDE} that are coherent with their associated Wiener paths. 


\paragraph{Standard splitting simulation procedure.} In the literature (see, e.g., \cite{Ableidinger2017,Buckwar2022}), realizations of $\eta_i(\Delta)$,  $i=1,\ldots,n$, required in the splitting schemes \eqref{eq:splitting}, are typically generated by sampling from the $d$-dimensional normal distribution $\mathcal{N}_d(0_d,C_\eta(\Delta))$ using the unique Cholesky factorization of the covariance matrix $C_\eta(\Delta)$. More precisely, one determines a matrix $L_\eta(\Delta)$ such that $L_\eta(\Delta) L_\eta(\Delta)^\top=C_\eta(\Delta)$, samples $\mathbf{z}$ from $\mathcal{N}_d(0_d,\mathbb{I}_{d\times d})$, and sets $\mathbf{x}=L_\eta(\Delta)\mathbf{z}$. Then $\mathbf{x}$ is a realization of $\eta_i(\Delta)$. However, as for linear SDEs (cf. Section \ref{sec2:LinearSDEs}) this approach does not preserve the coherence between the generated approximate path of SDE \eqref{eq:semi-linearSDE} and its underlying Wiener path. 


\paragraph{Wiener-path-based splitting simulation procedure.} To overcome this issue, and motivated by \cite{Cruz2020}, we propose the following procedure for sampling realizations of $\eta_i(\Delta)$, $i=1,\ldots,n$, required in~\eqref{eq:splitting}. Using the joint covariance matrix $C_{\xi;\eta}(\Delta)$ from \eqref{eq:C_xi_eta}, we determine a matrix $L$ satisfying $LL^\top=C_{\xi;\eta}(\Delta)$. Next, we sample $\mathbf{z}$ from $\mathcal{N}_{m+d}(0_{m+d},\mathbb{I}_{(m+d)\times(m+d)})$ and set $(\mathbf{y},\mathbf{x})^\top=L\mathbf{z}$. Then $(\mathbf{y},\mathbf{x})^\top$ is a joint realization of
\begin{equation*}
    \begin{pmatrix}
        \xi_i(\Delta) \\
        \eta_i(\Delta)
    \end{pmatrix},
\end{equation*}
where the $\xi_i(\Delta)$, $i=1,\ldots,n$, are as in \eqref{eq:xi_i}. 
The vector  $\mathbf{x}$ is then used in the splitting schemes \eqref{eq:splitting} to compute $\widetilde{X}(t_i)$ from $\widetilde{X}(t_{i-1})$, while $\mathbf{y}$ determines the Wiener increment via $W(t_i)=W(t_{i-1})+\mathbf{y}$. 

This procedure, which is an extension of Algorithm \ref{alg:Exact_WPB} for linear SDEs \eqref{eq:LinearSDE} to semi-linear SDEs~\eqref{eq:semi-linearSDE}, is summarized in  Algorithm \ref{alg:Strang_WPB} for the Strang splitting scheme. An analogous construction applies to the Lie-Trotter splitting scheme. Algorithm \ref{alg:Strang_WPB} returns both an approximate path of SDE \eqref{eq:semi-linearSDE} and its associated Wiener path. Note that, the latter could then be used to approximate the same path of SDE \eqref{eq:semi-linearSDE} with another numerical method, e.g., the Euler--Maruyama method. 

\begin{algorithm}[t]
	\caption{\textbf{Wiener-path-based} Strang splitting simulation (\textbf{WPB-Strang}) 
		\ \\ \textbf{Input:} Initial value $X_0$, time horizon $T$, 
  step-size $\Delta$, exponential matrix~$e^{A \Delta}$, joint covariance matrix $C_{\xi;\eta}(\Delta)$ and ODE solution function $f(\cdot;\cdot)$. 
  \ \\
		\textbf{Output:} Approximate path of the semi-linear SDE \eqref{eq:semi-linearSDE} and corresponding Wiener path. 
        }
	\label{alg:Strang_WPB}
    \begin{spacing}{1.1}
	\begin{algorithmic}[1]
        \State Set $n=T/\Delta$.
		\State Set $\widetilde{X}^{\textrm{S}}(t_0)=X_0$ and $W(t_0)=0_m$. 
        \State Determine $L$ such that $L L^\top=C_{\xi;\eta}(\Delta)$.
		\For{$i = 1:n$}
        \State Generate $\mathbf{z}=(z_1,\ldots,z_{m+d})^\top$ from $\mathcal{N}_{m+d}\bigl( 0_{m+d}, \mathbb{I}_{(m+d)\times (m+d)} \bigr)$.
        \State Set $(\mathbf{y},\mathbf{x})^\top=(y_1,\ldots,y_m,x_1,\ldots,x_d)^\top=L\mathbf{z}$.
        \State Set $\widetilde{X}^{\textrm{S}}(t_i)=f\Bigl(e^{A\Delta}f(\widetilde{X}^{\textrm{S}}(t_{i-1});\Delta/2)+\mathbf{x};\Delta/2\Bigr)$. 
        \State Set $W(t_i)=W(t_{i-1})+\mathbf{y}$.
		\EndFor
		\State Return $\widetilde{X}^{\textrm{S}}(t_i)$, 
        $i=0,\ldots,n$, and $W(t_i)$, $i=0,\ldots,n$. 
	\end{algorithmic}
    \end{spacing}
\end{algorithm}

\begin{remark}
    For the construction of the splitting schemes \eqref{eq:splitting}, we have used the exact solution $f$ of ODE \eqref{ODE}. Nevertheless, the Wiener-path-based simulation approach proposed in Algorithm \ref{alg:Strang_WPB} remains applicable when a numerical approximation of the solution to \eqref{ODE} is used instead.

    Moreover, note that the matrix exponential $e^{A\Delta}$ and the joint covariance matrix $C_{\xi;\eta}(\Delta)$ need to be computed only once for a given time step-size, thereby reducing the computational cost of the path simulation. However, deriving their closed-form expressions is generally nontrivial. In the subsequent section, we consider different stochastic models from neuroscience for which these expressions remain tractable and can be obtained in closed form.  
\end{remark}


\paragraph{Wiener-path-based splitting simulation on different time grids.}

\begin{algorithm}[t]
	\caption{\textbf{Wiener-path-based} Strang splitting on different time grids (\textbf{WPB-Strang-grids}) 
		\ \\ \textbf{Input:} Initial value $X_0$, time horizon $T$, 
  step-sizes $\Delta_1$ and $\Delta_2$, drift matrix $A$, exponential matrices $e^{A \Delta_1}$ and $e^{A \Delta_2}$, joint covariance matrix $C_{\xi;\eta}(\Delta_1)$ and ODE solution function $f(\cdot;\cdot)$. 
  \ \\
		\textbf{Output:} Approximate path of the semi-liner SDE \eqref{eq:semi-linearSDE} and corresponding Wiener path, generated on two different time grids. 
        }
	\label{alg:Strang_WPB_twoDeltas}
    \begin{spacing}{1.1}
	\begin{algorithmic}[1]
        \State Set $n_2=T/\Delta_2$, $l=\Delta_2/\Delta_1$ and $i=0$.
		\State Set $\widetilde{X}^{\textrm{S}}_{\Delta_1}(t_0)=\widetilde{X}^{\textrm{S}}_{\Delta_2}(t_0)=X_0$ and $W_{\Delta_1}(t_0)=W_{\Delta_2}(t_0)=0_m$. 
        \State Determine $L_{\Delta_1}$ such that $L_{\Delta_1} L_{\Delta_1}^\top=C_{\xi;\eta}(\Delta_1)$.
		\For{$j = 1:n_2$}
        \For{$k = 1:l$}
        \State Generate $\mathbf{z}_k$ from $\mathcal{N}_{m+d}\bigl( 0_{m+d}, \mathbb{I}_{(m+d)\times (m+d)} \bigr)$ and set $(\mathbf{y}_k,\mathbf{x}_k)^\top=L_{\Delta_1}\mathbf{z}_k$.
        \State Set $i=i+1$ (or $i=(j-1)l+k$).
		\State Set $\widetilde{X}^{\textrm{S}}_{\Delta_1}(t_i)=f\Bigl(e^{A\Delta_1}f(\widetilde{X}^{\textrm{S}}_{\Delta_1}(t_{i-1});\Delta_1/2)+\mathbf{x}_k;\Delta_1/2\Bigr)$. 
        \State Set $W_{\Delta_1}(t_i)=W_{\Delta_1}(t_{i-1})+\mathbf{y}_k$.
		\EndFor
        \State Set $\mathbf{y}=\sum\limits_{k=1}^{l}\mathbf{y}_k$.
        \State Set $\mathbf{x}=\sum\limits_{k=1}^{l} e^{A(l-k)\Delta_1} \mathbf{x}_k$.
        \State Set $\widetilde{X}^{\textrm{S}}_{\Delta_2}(t_j)=f\Bigl(e^{A\Delta_2}f(\widetilde{X}^{\textrm{S}}_{\Delta_2}(t_{j-1});\Delta_2/2)+\mathbf{x};\Delta_2/2\Bigr)$.
        \State Set $W_{\Delta_2}(t_j)=W_{\Delta_2}(t_{j-1})+\mathbf{y}$.
        \EndFor
		\State Return $\widetilde{X}^{\textrm{S}}_{\Delta_1}(t_i)$ and $W_{\Delta_1}(t_i)$, $i=0,\ldots,n_1$, and $\widetilde{X}^{\textrm{S}}_{\Delta_2}(t_j)$ and $W_{\Delta_2}(t_j)$, $j=0,\ldots,n_2$. 
	\end{algorithmic}
    \end{spacing}
\end{algorithm}

Finally, we aim to generate an approximate splitting realization of SDE \eqref{eq:semi-linearSDE} simultaneously on two different time grids, thereby extending the procedure introduced in Section \ref{sec2.2:LinearSDEs_WPB_exact_grids}. For the Euler--Maruyama method, this construction is straightforward. Indeed, the Wiener path $W_{\Delta_1}(t_i)$, $i=0,\ldots,n_1$, used to generate the approximate path on the fine grid can simply be subsampled to obtain the corresponding Wiener path $W_{\Delta_2}(t_j)$, $j=0,\ldots,n_2$, on the coarse grid. For the splitting schemes \eqref{eq:splitting}, however, the construction is less immediate, since both the approximate path of SDE \eqref{eq:semi-linearSDE} and its associated Wiener path are generated simultaneously by the Wiener-path-based simulation procedure, rather than by first prescribing a Wiener path and subsequently generating the corresponding approximate SDE path, as in the Euler--Maruyama scheme. Nevertheless, by again making use of Proposition~\ref{prop:coarse_eta_xi}, Algorithm \ref{alg:Exact_WPB_twoDeltas} can be  extended to the splitting setting, resulting in Algorithm \ref{alg:Strang_WPB_twoDeltas} for the Strang scheme in \eqref{eq:splitting}. An analogous construction applies again for the Lie-Trotter scheme.

\begin{remark}
    Note that Algorithm \ref{alg:Exact_WPB_twoDeltas} for linear SDEs \eqref{eq:LinearSDE} yields $X_{\Delta_2}(t_j)=X_{\Delta_1}(t_{lj})$, $j=1,\ldots,n_2$, as discussed in Section \ref{sec2.2:LinearSDEs_WPB_exact_grids}. In contrast, Algorithm \ref{alg:Strang_WPB_twoDeltas} for semi-linear SDEs \eqref{eq:semi-linearSDE} generally gives $\widetilde{X}^{\textrm{S}}_{\Delta_2}(t_j)\neq\widetilde{X}^{\textrm{S}}_{\Delta_1}(t_{lj})$, $j=1,\ldots,n_2$. This difference is expected, since the accuracy of the numerical method depends on the chosen time step $\Delta$. Nevertheless, Algorithm \ref{alg:Strang_WPB_twoDeltas} preserves the consistency of the Wiener paths across the two time grids, in the sense that $W_{\Delta_2}(t_j)=W_{\Delta_1}(t_{lj})$, $j=1,\ldots,n_2$, while each approximate SDE path remains associated with its corresponding underlying Wiener~path.
\end{remark}


\section{Illustration of the algorithms on stochastic models in neuroscience}
\label{sec5:Neuronal_models}

In this section, the proposed \textit{WPB-Strang} (Algorithm \ref{alg:Strang_WPB}) and \textit{WPB-Strang-grids} (Algorithm \ref{alg:Strang_WPB_twoDeltas}) algorithms are illustrated on two widely studied stochastic models in neuroscience: the stochastic FitzHugh--Nagumo (FHN) model and the stochastic Jansen--Rit neural mass model (JRNMM).


\subsection{Stochastic FitzHugh--Nagumo model}
\label{sec5.1:FHN_hypo}

The stochastic FHN model is a classical mathematical model for describing the dynamics of a single neuron, originating from the works of FitzHugh (1961) \cite{FitzHugh1961} and Nagumo et al. (1962)~\cite{Nagumo1962}. The model captures the generation of spikes, also referred to as action potentials, at the intracellular level under stochastic influences. In its \textit{elliptic} form (see, e.g., \cite{Berglund2012,Bonaccorsi2008,Buckwar2022}), it is given by the following $(d=2)$-dimensional non-linear SDE with additive noise:
\begin{equation}\label{eq:FHN_elliptic}
	d \underbrace{\begin{pmatrix}
			V(t) \\
			U(t) 
	\end{pmatrix}}_{=X(t)}
	=
	\underbrace{\begin{pmatrix}
			\frac{1}{\epsilon}\Bigl(V(t)-V^3(t)-U(t)\Bigr) \\
			\gamma V(t)-U(t)+\beta
	\end{pmatrix}}_{=F(X(t))} dt \ + \
	\underbrace{\begin{bmatrix}
		\sigma_1 &	0 \\ 
        0 & \sigma_2
	\end{bmatrix}}_{=\Sigma} dW(t), \qquad X(0)=X_0, 
\end{equation}
where $t\in[0,T]$, $W(t)=(W^1(t),W^2(t))^\top$, $t\in [0,T]$, is an $(m=2)$-dimensional standard Wiener process, and $\sigma_1,\sigma_2>0$ are corresponding noise intensity parameters. A commonly used \textit{hypoelliptic} formulation of the model (see, e.g., \cite{Ditlevsen2019,Leon2018,Samson2025}) is briefly discussed in Remark~\ref{rem:hypo_FHN} at the end of this subsection. 

The $V$-component of system \eqref{eq:FHN_elliptic} describes the evolution of the neuronal membrane voltage, which can be measured using suitable recording techniques, whereas the $U$-component is an unobserved recovery variable that models the ion channel kinetics. The stochastic FHN model~\eqref{eq:FHN_elliptic} contains five parameters: $\beta,\gamma>0$, which govern the position and duration of an action potential, respectively; $\epsilon>0$, which determines the time-scale separation of the system; and $\sigma_1,\sigma_2>0$, which control the noise intensity.


\paragraph{Wiener-path-based Strang splitting for the stochastic FHN model.}

Following \cite{Buckwar2022}, we decompose the drift coefficient $F(X(t))$ of the stochastic FHN model \eqref{eq:FHN_elliptic} into a linear and a non-linear part according to  $F(X(t))=AX(t)+N(X(t))$, where 
\begin{equation}\label{eq:splitting_choice_FHN}
    A=\begin{bmatrix}
        0 & -\frac{1}{\epsilon} \\
        \gamma & -1
    \end{bmatrix}, \qquad N(X(t))=\begin{pmatrix}
        \frac{1}{\epsilon}(V(t)-V^3(t)) \\
        \beta
    \end{pmatrix}.
\end{equation}
These terms define the two subequations \eqref{SDE} and \eqref{ODE} within the splitting framework.

\begin{remark}
    Note that the matrix $A$ in \eqref{eq:splitting_choice_FHN} is chosen such that subequation \eqref{SDE} is hypoelliptic (cf. Remark~\ref{rem:A_assumptions}). Moreover, the function $N$ in \eqref{eq:splitting_choice_FHN} is globally one-sided Lipschitz continuous and satisfies a polynomial growth condition (cf. Remark \ref{rem:N_assumptions}). The interested reader is referred to \cite{Buckwar2022} for further details in this regard.
\end{remark}

In the following, we provide closed-form expressions for the solution $f(\cdot;\cdot)$ of the deterministic subequation \eqref{ODE}, and the matrix exponential $e^{A\Delta}$ and covariance matrix $C_\eta(\Delta)$, which appear in the solution of the stochastic subequation \eqref{SDE}. In addition, we derive the joint covariance matrix $C_{\xi,\eta}(\Delta)$, which contains $C_\eta(\Delta)$ as a submatrix and enables  Wiener-path-based simulation via Algorithms \ref{alg:Strang_WPB} and \ref{alg:Strang_WPB_twoDeltas}, respectively.

For $t>0$ and $x=(v,u)^\top \in \mathbb{R}^2$, the function $f$ is given by 
\begin{equation}\label{eq:f_FHN}
     f(x;t)=\begin{pmatrix}
		v/{\sqrt{e^{-\frac{2t}{\epsilon}}+v^2\left(1-e^{-\frac{2t}{\epsilon}}\right)}} \\
		 \beta t + u
	\end{pmatrix}.
\end{equation}

Next, define   
\begin{equation}\label{eq:kappa}
    \kappa:= \frac{4\gamma}{\epsilon}-1.
\end{equation}
Then, for $t>0$ and $\kappa>0$, the matrix exponential is given by
\begin{equation}\label{eq:MatrixExp_FHN}
	e^{A t}=e^{-\frac{t}{2}}\begin{bmatrix}
		\cos(\frac{1}{2}\sqrt{\kappa}t)+\frac{1}{\sqrt{\kappa}}\sin(\frac{1}{2}\sqrt{\kappa}t) & -\frac{2}{\epsilon\sqrt{\kappa}}\sin(\frac{1}{2}\sqrt{\kappa}t) \\[0.5em]
		\frac{2\gamma}{\sqrt{\kappa}}\sin(\frac{1}{2}\sqrt{\kappa}t) & \cos(\frac{1}{2}\sqrt{\kappa}t)-\frac{1}{\sqrt{\kappa}}\sin(\frac{1}{2}\sqrt{\kappa}t)
	\end{bmatrix}.
\end{equation}
Moreover, for $t>0$ and $\kappa>0$, the joint covariance matrix takes the form
\begin{equation}\label{eq:JointCovMat_FHN}
    C_{\xi;\eta}(t)=\begin{bmatrix}
        C_\xi(t) & C(t)^\top \\
        C(t) & C_\eta(t)
    \end{bmatrix},
\end{equation}
where
\begin{align*}
    C_\xi(t)&=t\mathbb{I}_{2\times 2}, \\[1em]
    C_\eta(t)
    &=\begin{bmatrix}
        c_{\eta;11}(t) & c_{\eta;12(t)} \\
        c_{\eta;12(t)} & c_{\eta;22(t)}
    \end{bmatrix}, \quad \text{with} \\[1em]
    & c_{\eta;11}(t) = \frac{\epsilon e^{-t}}{2\gamma\kappa} \biggl( -\frac{4\gamma}{\epsilon^2}( \sigma_1^2\gamma+\sigma_2^2\frac{1}{\epsilon})  + \kappa e^t \Bigl(\sigma_1^2(1+\frac{\gamma}{\epsilon}) +\sigma_2^2\frac{1}{\epsilon^2}\Bigr) \\ \nonumber & \hspace{1.5cm} + \Bigl( \sigma_1^2(1-\frac{3\gamma}{\epsilon}) +\sigma_2^2\frac{1}{\epsilon^2} \Bigr)\cos(\sqrt{\kappa}t) - \sqrt{\kappa} \Bigl(\sigma_1^2(1-\frac{\gamma }{\epsilon})+\sigma_2^2\frac{1}{\epsilon^2}\Bigr) \sin(\sqrt{\kappa}t) \biggr), \\[1em]
    & c_{\eta;12}(t) = \frac{\epsilon e^{-t}}{2 \kappa } \biggl( \sigma_1^2\kappa e^t - \frac{2}{\epsilon}(\sigma_1^2 \gamma  +\sigma_2^2\frac{1}{\epsilon}) \\ \nonumber & \hspace{1.5cm} +
    \Bigl( \sigma_1^2(1-\frac{2\gamma}{\epsilon}) + 2\sigma_2^2 \frac{1}{\epsilon^2} \Bigr) \cos(\sqrt{\kappa}t)  -\sigma_1^2\sqrt{\kappa} \sin(\sqrt{\kappa}t)  \biggr), \\[1em]
    & c_{\eta;22}(t) = \frac{\epsilon e^{-t}}{2\kappa} \biggl( (\sigma_2^2\frac{1}{\epsilon}+\sigma_1^2\gamma) \Bigl(  \cos(\sqrt{\kappa}t)-\frac{4\gamma}{\epsilon}+\kappa e^t \Bigr) + (\sigma_2^2\frac{1}{\epsilon}-\sigma_1^2\gamma) \sqrt{\kappa} \sin(\sqrt{\kappa}t) \biggr), \\[1em]
    C(t)
    &=\begin{bmatrix}
        c_{11}(t) & c_{12(t)} \\
        c_{21(t)} & c_{22(t)}
    \end{bmatrix}, \quad \text{with} \\[1em]
    & c_{11}(t) = \frac{\sigma_1\epsilon}{\gamma} \left( 1+e^{-t/2} \Bigl( -\cos( \frac{1}{2}\sqrt{\kappa}t ) + \frac{1}{\sqrt{\kappa}} (\frac{2\gamma}{\epsilon}-1) \sin( \frac{1}{2}\sqrt{\kappa}t)  \Bigr)  \right), \\[1em]
    & c_{12}(t) = \frac{\sigma_2}{\gamma} \left( -1+e^{-t/2} \Bigl( \cos( \frac{1}{2}\sqrt{\kappa}t ) + \frac{1}{\sqrt{\kappa}}\sin( \frac{1}{2}\sqrt{\kappa}t)  \Bigr)  \right), \\[1em]
    & c_{21}(t) = \sigma_1\epsilon \left( 1+e^{-t/2} \Bigl( -\cos( \frac{1}{2}\sqrt{\kappa}t ) - \frac{1}{\sqrt{\kappa}}\sin( \frac{1}{2}\sqrt{\kappa}t)  \Bigr)  \right), \\[1em]
    & c_{22}(t) = \frac{2\sigma_2 e^{-t/2}}{\sqrt{\kappa}}\sin\left(\frac{1}{2}\sqrt{\kappa}t \right). 
\end{align*}

\begin{remark}
    The expressions \eqref{eq:MatrixExp_FHN} and \eqref{eq:JointCovMat_FHN} are stated for the case $\kappa>0$, where $\kappa$ is defined in \eqref{eq:kappa}. This condition ensures that subequation \eqref{SDE} corresponds to a weakly damped stochastic harmonic oscillator, thereby producing oscillatory neural spiking dynamics. The interested reader is referred to \cite{Buckwar2022} for further details and the corresponding expressions in the cases $\kappa<0$ and $\kappa=0$.
\end{remark}


\paragraph{Numerical experiments for the stochastic FHN model.} 

Algorithm \ref{alg:Strang_WPB} is now applied to the stochastic FHN model \eqref{eq:FHN_elliptic} using the input components $f(\cdot;\cdot)$ from \eqref{eq:f_FHN}, $e^{A \Delta}$ from \eqref{eq:MatrixExp_FHN}, and $C_{\xi;\eta}(\Delta)$ from  \eqref{eq:JointCovMat_FHN}. The initial value is chosen as $X_0=(-1.1,1)^\top$, the time horizon is $T=5$, and the step-size is $\Delta=10^{-4}$. The model parameters are set to $\beta=0.1$, $\gamma=5$, $\epsilon=0.05$, and $\sigma_1=\sigma_2=0.5$.

\begin{figure}[t]
	\begin{centering}
		\includegraphics[width=1.0\textwidth]{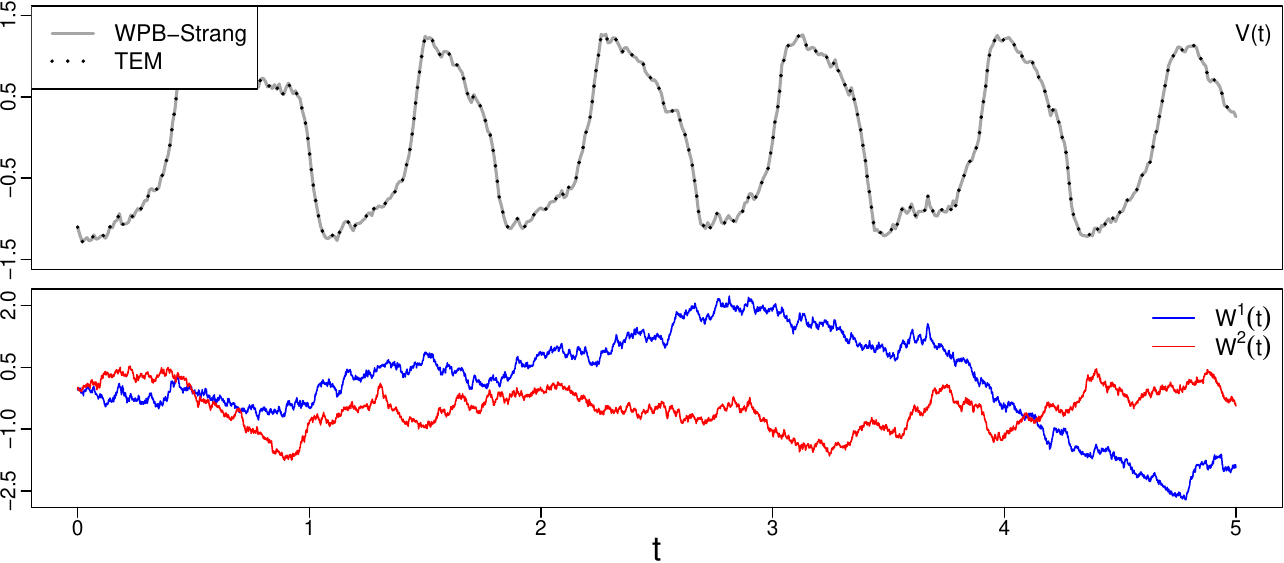}	
		\caption{Path of the $V$-component of the stochastic FHN model \eqref{eq:FHN_elliptic}. Top panel: realization $\widetilde{V}^{\textrm{S}}(t_i)$ (gray solid line) obtained via Algorithm \ref{alg:Strang_WPB}, compared with the tamed Euler--Maruyama realization $\widetilde{V}^{\textrm{TEM}}(t_i)$ (black dotted line) obtained using the same underlying Wiener path. Bottom panel: corresponding Wiener path $W(t_i)=(W^1(t_i),W^2(t_i))^\top$, $i=0,\ldots,n$.}
		\label{fig:FHN_elliptic}
	\end{centering}
\end{figure}

Figure \ref{fig:FHN_elliptic} shows the resulting path $\widetilde{V}^{\textrm{S}}(t_i)$, $i=0,\ldots,n$, of the $V$-component (top panel; gray solid line), together with the corresponding underlying Wiener path $W(t_i)=(W^1(t_i),W^2(t_i))^\top$, $i=0,\ldots,n$ (bottom panel; blue and red solid lines).  

In addition, the Wiener path $W(t_i)$, $i=0,\ldots,n$, returned by Algorithm \ref{alg:Strang_WPB}, is used to generate a corresponding realization of the $V$-component via the tamed Euler--Maruyama (TEM) method~\cite{Hutzenthaler2012_2}, which incorporates a taming term to handle the superlinearly growing drift. 
The resulting path $\widetilde{V}^{\textrm{TEM}}(t_i)$, $i=0,\ldots,n$, obtained from the TEM scheme using the Wiener path $W(t_i)$, $i=0,\ldots,n$, generated by Algorithm \ref{alg:Strang_WPB}, is also shown in the top panel of Figure \ref{fig:FHN_elliptic} (black dotted line). As expected, the path overlaps with $\widetilde{V}^{\textrm{S}}(t_i)$, $i=0,\ldots,n$, produced by Algorithm \ref{alg:Strang_WPB}, due to the small time step $\Delta=10^{-4}$ used in both simulations.

\begin{figure}[t]
	\begin{centering}
		\includegraphics[width=1.0\textwidth]{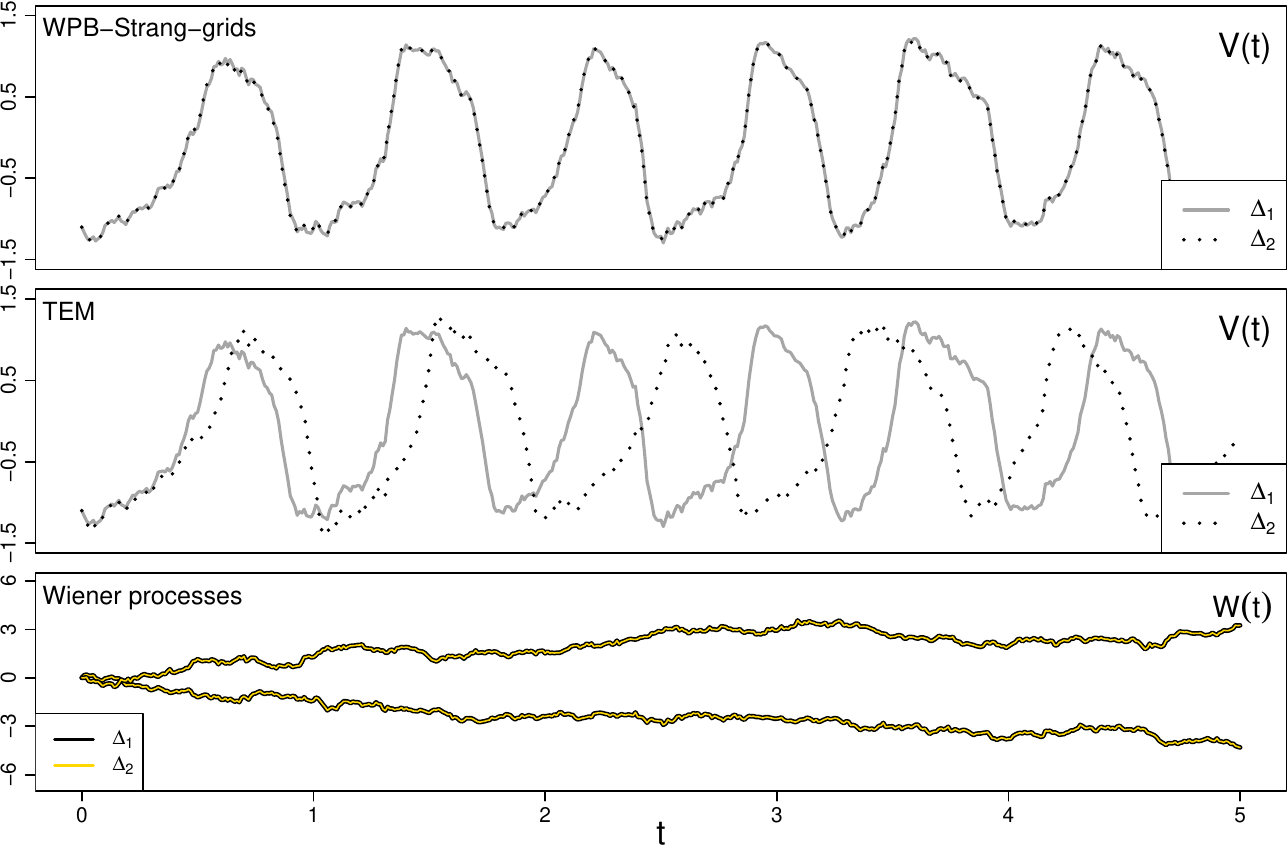}	
		\caption{Path of the $V$-component of the stochastic FHN model \eqref{eq:FHN_elliptic}. Top panel: fine- and coarse-grid paths  $\widetilde{V}^{\textrm{S}}_{\Delta_1}(t_i)$ (gray solid line) and $\widetilde{V}^{\textrm{S}}_{\Delta_2}(t_j)$ (black dotted line) obtained via Algorithm~\ref{alg:Strang_WPB_twoDeltas}. Middle panel: corresponding tamed Euler--Maruyama realizations $\widetilde{V}^{\textrm{TEM}}_{\Delta_1}(t_i)$ (gray solid line) and $\widetilde{V}^{\textrm{TEM}}_{\Delta_2}(t_j)$ (black dotted line). Bottom panel: underlying Wiener paths $W_{\Delta_1}(t_i)$ (black lines) and $W_{\Delta_2}(t_j)$ (yellow lines).}
		\label{fig:FHN_elliptic_twoDeltas}
	\end{centering}
\end{figure}

We now apply Algorithm \ref{alg:Strang_WPB_twoDeltas} to the stochastic FHN model \eqref{eq:FHN_elliptic}, using the same algorithm input parameters and components as before as well as the step-sizes $\Delta_1=10^{-4}$ and $\Delta_2=10^{-2}$. 

Figure \ref{fig:FHN_elliptic_twoDeltas} shows the resulting paths $\widetilde{V}^{\textrm{S}}_{\Delta_1}(t_i)$, $i=0,\ldots,n_1$, and $\widetilde{V}^{\textrm{S}}_{\Delta_2}(t_j)$, $j=0,\ldots,n_2$, on the fine and coarse grids, respectively (top panel: gray solid and black dotted lines). The associated Wiener paths
\begin{equation*}
    W_{\Delta_1}(t_i)=(W^1_{\Delta_1}(t_i),W^2_{\Delta_1}(t_i))^\top, \quad i=0,\ldots,n_1,
\end{equation*}
and
\begin{equation*}
    W_{\Delta_2}(t_j)=(W^1_{\Delta_2}(t_j),W^2_{\Delta_2}(t_j))^\top, \quad j=0,\ldots,n_2,
\end{equation*}
are displayed in the bottom panel (black and yellow lines, respectively). 

In addition, the Wiener paths $W_{\Delta_1}(t_i)$ and $W_{\Delta_2}(t_j)$, returned by Algorithm \ref{alg:Strang_WPB_twoDeltas}, are used to generate the corresponding realizations $\widetilde{V}^{\textrm{TEM}}_{\Delta_1}(t_i)$ and $\widetilde{V}^{\textrm{TEM}}_{\Delta_2}(t_j)$ via the TEM scheme. These paths are shown in the middle panel of Figure \ref{fig:FHN_elliptic_twoDeltas} (gray solid and black dotted line, respectively).

It can be observed that the splitting scheme produces a coarse-grid path for $\Delta_2=10^{-2}$ that closely overlaps with the corresponding fine-grid path for $\Delta_1=10^{-4}$. In contrast, this behaviour is not observed for the TEM scheme, where the path corresponding to $\Delta_2=10^{-2}$ deviates visibly from that obtained with $\Delta_1=10^{-4}$. This observation is consistent with previous findings in the literature, where splitting methods for the stochastic FHN model have been reported to outperform the TEM scheme (and related EM-type methods) in terms of structure preservation, including the preservation of amplitudes, frequencies, and phases of oscillations; see, e.g., \cite{Buckwar2022,Fosteretal2022}. The results thus highlight the advantage of Strang splitting, which allows for larger time steps while maintaining close agreement with the fine-grid path.

\begin{remark}
\label{rem:hypo_FHN}
    The \textit{hypoelliptic} stochastic FHN model is obtained from \eqref{eq:FHN_elliptic} by setting $\sigma_1=0$, such that the system is driven by an $(m=1)$-dimensional Wiener process. The drift decomposition~\eqref{eq:splitting_choice_FHN} remains unchanged, and hence the expressions for $f(\cdot;\cdot)$ in \eqref{eq:f_FHN} and $e^{At}$ in \eqref{eq:MatrixExp_FHN} remain valid. Only the joint covariance matrix \eqref{eq:JointCovMat_FHN} needs to be adapted: $C_\xi(t)=t$, $C_\eta(t)$ is obtained by setting $\sigma_1=0$ in the expressions above, and $C(t)$ reduces to its second column, yielding a $(3\times 3)$-dimensional joint covariance matrix $C_{\xi;\eta}(t)$. With these modifications, Algorithms \ref{alg:Strang_WPB} and \ref{alg:Strang_WPB_twoDeltas} apply analogously to the hypoelliptic case, producing closely overlapping Wiener-path-based realizations also in this setting (results not shown).
\end{remark}


\subsection{Stochastic Jansen--Rit neural mass model}
\label{sec5.3:JRNMM}

The stochastic JRNMM originates from the work of Jansen and Rit (1995) \cite{Jansen1995}; see also, e.g., \cite{Ableidinger2017,Ditlevsenetal2025,Wendling2000} for further extensions and discussions. The model describes the electrical activity of a population of neurons in terms of average population dynamics and captures the interaction between main pyramidal cells and excitatory and inhibitory interneurons. It has been shown to successfully reproduce electroencephalogram (EEG) signals. The model is given by the following ($d=6$)-dimensional non-linear SDE with additive noise:
\begin{equation}\label{eq:JRNMM}
\begin{aligned}
	dX_1(t)&=X_4(t)dt,\\
	dX_2(t)&=X_5(t)dt,\\
	dX_3(t)&=X_6(t)dt,\\
	dX_4(t)&=\left[Aa\bigl(\textrm{sig}\left(X_2(t)-X_3(t)\right) \bigr) -2aX_4(t)-a^2X_1(t)\right]  dt+ \sigma_1 dW^1(t),\\
	dX_5(t)&=\left[Aa \bigl(\mu+C_2\textrm{sig}\bigl(C_1X_1(t)\bigr)\bigr)-2aX_5(t)-a^2X_2(t)\right]dt+\sigma_2 dW^2(t),\\
	dX_6(t)&= \left[BbC_4\textrm{sig}\left(C_3 X_1(t)\right)-2bX_6(t)-b^2X_3(t)\right]dt+\sigma_3 dW^3(t),
\end{aligned}
\end{equation}
where $t \in [0,T]$, $W(t)=(W^1(t),W^2(t),W^3(t))^\top$, $t\in [0,T]$, is an ($m=3$)-dimensional standard Wiener process, and the sigmoid function $\textrm{sig} \colon \mathbb{R} \to [0,v_{\textrm{max}}]$ is defined by
\begin{equation}\label{eq:sigmoid}
    \textrm{sig}(x):=\frac{v_{\textrm{max}}}{1+e^{r(v_0-x)}}.
\end{equation}

The model component 
\begin{equation*}
    X_2(t)-X_3(t), \qquad t \in [0,T],
\end{equation*}
corresponds to the observable EEG signal, whereas the individual state variables $X_1(t),\ldots,X_6(t)$, $t\in [0,T]$, are not directly observable. The stochastic JRNMM \eqref{eq:JRNMM} contains several biologically motivated parameters ($A$, $a$, $B$, $b$, $C_1$, $C_2$, $C_3$, $C_4$, $\mu$, $\sigma_1$, $\sigma_2$, $\sigma_3$, $v_{\textrm{max}}$, $v_0$, $r$) whose interpretations and typical values are well documented in the literature; see, for example, \cite{Ableidinger2017}, and in particular Table $1$ therein.

\paragraph{Wiener-path-based Strang splitting for the stochastic JRNMM.}

The splitting method for the stochastic JRNMM \eqref{eq:JRNMM} is based on a reformulation of the model as a stochastic Hamiltonian-type system; see \cite{Ableidinger2017}. To this end, define $X(t):=(Q(t),P(t))^\top$, where $Q(t)=(X_1(t),X_2(t),X_3(t))^\top$ and $P(t)=(X_4(t),X_5(t),X_6(t))^\top$. Furthermore, introduce the $3\times 3$-dimensional diagonal matrices  
$\Gamma:=\textrm{diag}[a,a,b]$ and $\Lambda:=\textrm{diag}[\sigma_1,\sigma_2,\sigma_3]$. In addition, define the function $G\colon \mathbb{R}^3 \to \mathbb{R}^3$ by
\begin{equation*}
    G(Q(t)):=\begin{pmatrix}
        Aa\bigl(\textrm{sig}\left(X_2(t)-X_3(t)\right) \bigr) \\
        Aa \bigl(\mu+C_2\textrm{sig}\bigl(C_1X_1(t)\bigr)\bigr) \\
        BbC_4\textrm{sig}\left(C_3 X_1(t)\right)
    \end{pmatrix}.
\end{equation*}
With these definitions, the stochastic JRNMM \eqref{eq:JRNMM} can be written as a Hamiltonian-type SDE of the form
\begin{equation}\label{eq:JRNMM_Hamiltonian}
        d\underbrace{\begin{pmatrix}
            Q(t) \\
            P(t)
        \end{pmatrix}}_{=X(t)}=\underbrace{\begin{pmatrix}
            P(t) \\
            -\Gamma^2Q(t)-2\Gamma P(t) + G(Q(t))
        \end{pmatrix}}_{=F(X(t))} dt + \underbrace{\begin{bmatrix}
            \mathbb{O}_{3\times 3} \\ 
            \Lambda
        \end{bmatrix}}_{=\Sigma} dW(t), \qquad X(0)=X_0.
\end{equation}

Following \cite{Ableidinger2017}, we now decompose the drift coefficient $F(X(t))$ of the Hamiltonian-type system~\eqref{eq:JRNMM_Hamiltonian} into a linear and non-linear part according to $F(X(t))=AX(t)+N(X(t))$, where
\begin{equation}\label{eq:splitting_choice_JRNMM}
    A=\begin{bmatrix}
        \mathbb{O}_{3\times 3} & \mathbb{I}_{3\times 3} \\
        -\Gamma^2 & -2 \Gamma
    \end{bmatrix}, \qquad N(X(t))= \begin{pmatrix}
        0_3 \\
        G(Q(t))
    \end{pmatrix}.
\end{equation}
These terms define the two subequations \eqref{ODE} and \eqref{SDE} within the splitting framework.

\begin{remark}
    The matrix $A$ in \eqref{eq:splitting_choice_JRNMM} is again chosen such that subequation \eqref{SDE} is hypoelliptic (cf.~Remark~\ref{rem:A_assumptions}). Moreover, since the sigmoid function defined in \eqref{eq:sigmoid} is globally Lipschitz continuous, the function $N$ in \eqref{eq:splitting_choice_JRNMM} is globally Lipschitz continuous as well (cf. Remark \ref{rem:N_assumptions}). The reader is referred to  \cite{Ableidinger2017} for further details.
\end{remark}

In the following, closed-form expressions are provided for the solution $f(\cdot;\cdot)$ of the deterministic subequation \eqref{ODE}, and the matrix exponential $e^{A\Delta}$ and covariance matrix $C_\eta(\Delta)$, which appear in the solution of the stochastic subequation \eqref{SDE}. In addition, the joint covariance matrix $C_{\xi;\eta}(\Delta)$ is derived, which contains $C_\eta(\Delta)$ as a submatrix and enables Wiener-path-based simulation of the stochastic JRNMM via Algorithms \ref{alg:Strang_WPB} and \ref{alg:Strang_WPB_twoDeltas}, respectively.

For $t>0$ and $x=(q,p)^\top \in \mathbb{R}^6$, the function $f$ is given by
\begin{equation}\label{eq:f_JRNMM}
    f(x;t)= x+t \begin{pmatrix}
        0_3 \\
        G(q)
    \end{pmatrix}.
\end{equation}

Moreover, for $t>0$, the matrix exponential is given by
\begin{equation}\label{eq:MatrixExp_JRNMM}
    e^{At}= \begin{bmatrix}
        e^{-\Gamma t}\left(\mathbb{I}_{3\times 3}+\Gamma t\right) & e^{-\Gamma t} t\\
        -\Gamma^2 e^{-\Gamma t} t & e^{-\Gamma t}\left(\mathbb{I}_{3\times 3}-\Gamma t \right)
    \end{bmatrix}=: \begin{bmatrix}
        \nu(t) & \zeta(t) \\
        \nu^{\prime}(t) & \zeta^{\prime}(t)
    \end{bmatrix},
\end{equation}
while the joint covariance matrix takes the form
\begin{equation}\label{eq:JointCovMat_JRNMM}
    C_{\xi;\eta}(t)=\begin{bmatrix}
        C_\xi(t) & C(t)^\top \\
        C(t) & C_\eta(t)
    \end{bmatrix},
\end{equation}
where
\begin{align*}
    C_\xi(t)&=t \mathbb{I}_{3\times 3}, \\[1em]
    C_\eta(t)
    &= \begin{bmatrix}
        \frac{1}{4}\Gamma^{-3} \Lambda^2 \left( \mathbb{I}_{3\times 3} +\zeta(t)\nu^{\prime}(t) - \nu(t)^2 \right) & \frac{1}{2}\Lambda^2 \zeta(t)^2 \\
        \frac{1}{2}\Lambda^2 \zeta(t)^2 & \frac{1}{4}\Gamma^{-1} \Lambda^2 \left( \mathbb{I}_{3\times 3} +\zeta(t)\nu^{\prime}(t) - \zeta^{\prime}(t)^2 \right)
    \end{bmatrix}  ,\\[1em]
    C(t)
    &= \begin{bmatrix}
        \Gamma^{-2} \Lambda \left( \mathbb{I}_{3\times 3}-\nu(t) \right) \\
        \Lambda \zeta(t)
    \end{bmatrix}
    . 
\end{align*}

\paragraph{Numerical experiments for the stochastic JRNMM.} Algorithm \ref{alg:Strang_WPB} is now applied to the stochastic JRNMM \eqref{eq:JRNMM} using the input components $f(\cdot;\cdot)$ from \eqref{eq:f_JRNMM}, $e^{A\Delta}$ from \eqref{eq:MatrixExp_JRNMM}, and $C_{\xi;\eta}(\Delta)$ from \eqref{eq:JointCovMat_JRNMM}. The initial value is chosen as $X_0=(0.08,18,15,-0.5,0,0)^\top$, the time horizon is $T=3$, and the step-size is $\Delta=10^{-5}$. The model parameters are set according to the typical values reported in the literature (see, e.g., Table $1$ in \cite{Ableidinger2017}), and $\mu=220$, $\sigma_1=\sigma_3=10$ and $\sigma_2=2000$. 

\begin{figure}[htbp]
	\begin{centering}
		\includegraphics[width=1.0\textwidth]{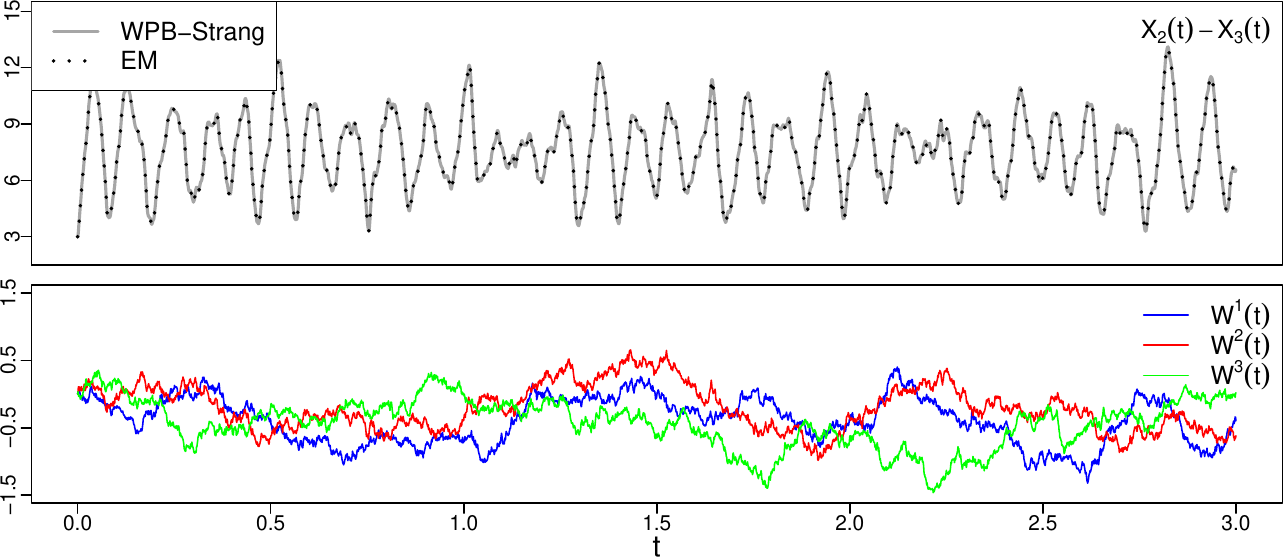}	
		\caption{Path of the $(X_2-X_3)$-component of the stochastic JRNMM \eqref{eq:JRNMM}. Top panel: realization $\widetilde{X}^{\textrm{S}}_2(t_i)-\widetilde{X}^{\textrm{S}}_3(t_i)$ (gray solid line) obtained via Algorithm \ref{alg:Strang_WPB}, compared with the Euler--Maruyama realization \mbox{$\widetilde{X}_2^{\textrm{EM}}(t_i)-\widetilde{X}_3^{\textrm{EM}}(t_i)$} (black dotted line) obtained using the same underlying Wiener path. Bottom panel: corresponding Wiener path $W(t_i)=(W^1(t_i),W^2(t_i),W^3(t_i))^\top$, $i=0,\ldots,n$.}
		\label{fig:JRNMM}
	\end{centering}
\end{figure}

\begin{figure}[htbp]
	\begin{centering}
		\includegraphics[width=1.0\textwidth]{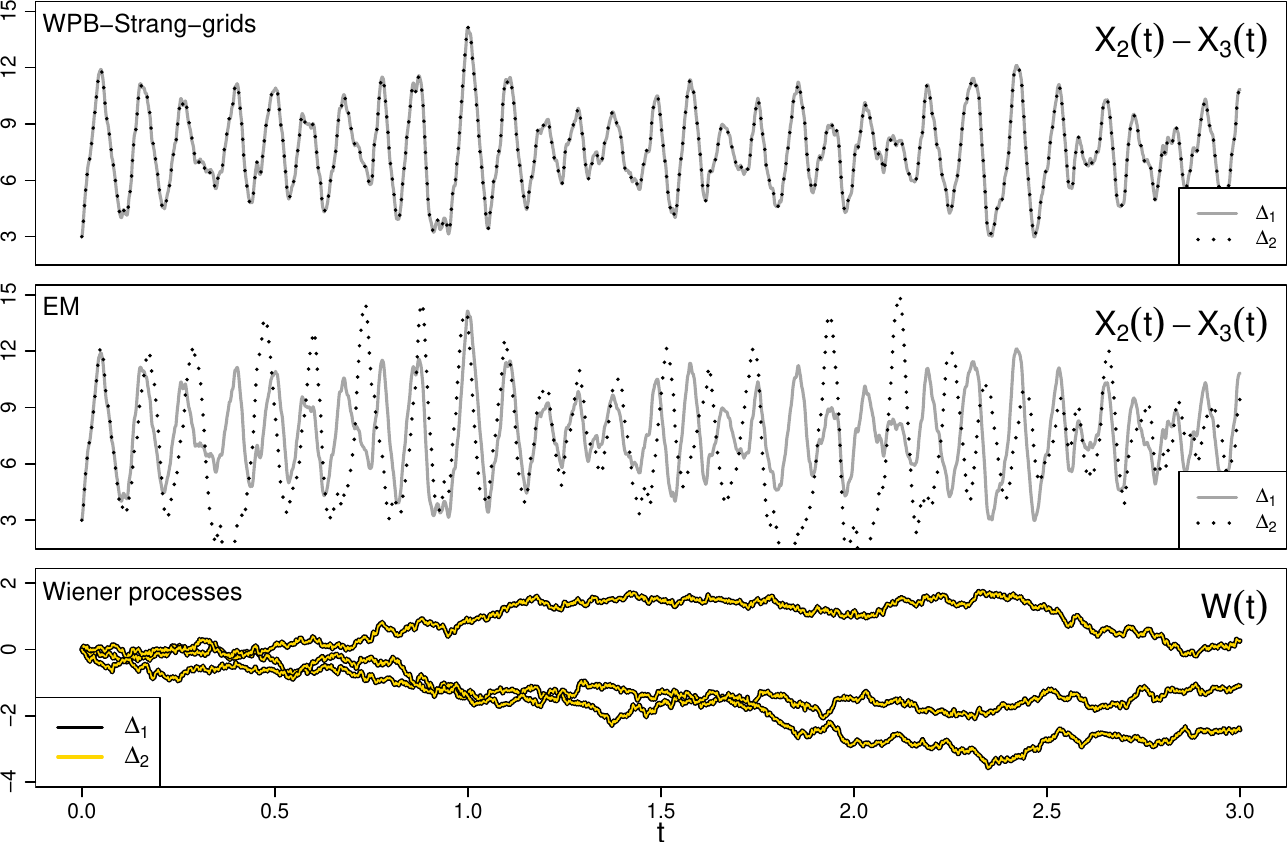}	
		\caption{Path of the $(X_2-X_3)$-component of the stochastic JRNMM \eqref{eq:JRNMM}. Top panel: fine- and coarse-grid paths  $\widetilde{X}^{\textrm{S}}_{2;\Delta_1}(t_i)-\widetilde{X}^{\textrm{S}}_{3;\Delta_1}(t_i)$ (gray solid line) and $\widetilde{X}^{\textrm{S}}_{2;\Delta_2}(t_j)-\widetilde{X}^{\textrm{S}}_{3;\Delta_2}(t_j)$ (black dotted line) obtained via Algorithm~\ref{alg:Strang_WPB_twoDeltas}. Middle panel: corresponding Euler--Maruyama realizations. Bottom panel: underlying Wiener paths $W_{\Delta_1}(t_i)$ (black lines) and $W_{\Delta_2}(t_j)$ (yellow lines).}
		\label{fig:JRNMM_twoDeltas}
	\end{centering}
\end{figure}

Figure \ref{fig:JRNMM} shows the resulting path $\widetilde{X}^{\textrm{S}}_2(t_i)-\widetilde{X}^{\textrm{S}}_3(t_i)$, $i=0,\ldots,n$, which describes alpha rhythmic EEG activity (top panel; gray solid line), together with the corresponding underlying Wiener path $W(t_i)=(W^1(t_i),W^2(t_i),W^3(t_i))^\top$, $i=0,\ldots,n$ (bottom panel; blue, red and green solid lines).

In addition, the Wiener path $W(t_i)$, $i=0,\ldots,n$, returned by Algorithm \ref{alg:Strang_WPB}, is used to generate a corresponding realization with the standard Euler--Maruyama (EM) method. The resulting path $\widetilde{X}_2^{\textrm{EM}}(t_i)-\widetilde{X}_3^{\textrm{EM}}(t_i)$, $i=0,\ldots,n$, obtained from the EM scheme using the Wiener path $W(t_i)$, $i=0,\ldots,n$, returned by Algorithm~\ref{alg:Strang_WPB}, is also shown in the top panel of Figure \ref{fig:JRNMM} (black dotted line). As expected, it overlaps with the path $\widetilde{X}^{\textrm{S}}_2(t_i)-\widetilde{X}^{\textrm{S}}_3(t_i)$, $i=0,\ldots,n$, generated by Algorithm~\ref{alg:Strang_WPB}, due to the small time step $\Delta=10^{-5}$ used in both simulations.

Algorithm \ref{alg:Strang_WPB_twoDeltas} is also applied to the stochastic JRNMM \eqref{eq:JRNMM}. The input paramters and components are the same as before, and the step-sizes are $\Delta_1=10^{-5}$ and $\Delta_2=10^{-3}$. 

Figure \ref{fig:JRNMM_twoDeltas} shows the resulting realizations 
\begin{equation*}
    \widetilde{X}^{\textrm{S}}_{2;\Delta_1}(t_i)-\widetilde{X}^{\textrm{S}}_{3;\Delta_1}(t_i), \ i=0,\ldots,n_1, \quad \text{and} \quad \widetilde{X}^{\textrm{S}}_{2;\Delta_2}(t_j)-\widetilde{X}^{\textrm{S}}_{3;\Delta_2}(t_j), \ j=0,\ldots,n_2,
\end{equation*}
(top panel: gray solid and black dotted lines, respectively), together with the corresponding Wiener paths 
\begin{equation*}
    W_{\Delta_1}(t_i)=(W^1_{\Delta_1}(t_i),W^2_{\Delta_1}(t_i),W^3_{\Delta_1}(t_i))^\top, \quad i=0,\ldots,n_1, 
\end{equation*}
and
\begin{equation*}
   W_{\Delta_2}(t_j)=(W^1_{\Delta_2}(t_j),W^2_{\Delta_2}(t_j),W^3_{\Delta_2}(t_j))^\top, \quad j=1,\ldots,n_2,
\end{equation*}
(bottom panel: black and yellow lines, respectively).

It can be observed that the approximate paths produced by the splitting method overlap closely, whereas the corresponding paths obtained via the EM scheme (middle panel: gray solid and black dotted lines, respectively) deviate significantly from one another. As for the stochastic FHN model, this observation agrees with results reported in the literature, where splitting methods for the stochastic JRNMM have been shown to outperform the EM scheme in preserving structural properties of the model; see, e.g., \cite{Ableidinger2017,Buckwar2019,Ditlevsenetal2025}. These results not only further highlight the advantage of splitting methods in allowing for larger time steps, but also the relevance of the proposed Algorithms \ref{alg:Strang_WPB} and \ref{alg:Strang_WPB_twoDeltas}, which enable a direct comparison of paths produced with different numerical methods and time steps.


\section{An alternative simulation procedure for the case $d=m$}
\label{appA}

In this section, we discuss an alternative to Algorithm \ref{alg:Exact_WPB} for the case $d=m$ of the linear SDE~\eqref{eq:LinearSDE}, which constructs the required random realizations from a single $d$-dimensional normal sample instead of a sample from a joint $(m+d)$-dimensional normal distribution. In Section \ref{sec6:LinSDEs}, we present the algorithm along with the corresponding theoretical results and illustrate them on the Ornstein--Uhlenbeck (OU) process and the simple SHO (cf. Section \ref{sec3:simple_SHO}). In Section \ref{sec6:SemiLinSDEs}, we then extend this approach to semi-linear SDEs of the form \eqref{eq:semi-linearSDE} and illustrate it on the elliptic stochastic FHN~model (cf. Section \ref{sec5.1:FHN_hypo}). 

\subsection{Linear SDEs with additive noise for $d=m$}\label{sec6:LinSDEs}

Algorithm \ref{alg:Exact_WPB} for linear SDEs of the form \eqref{eq:LinearSDE} preserves the joint covariance structure $C_{\xi;\eta}(\Delta)$ from~\eqref{eq:C_xi_eta} by sampling realizations~$\textbf{z}$ from $\mathcal{N}_{m+d}(0_{m+d},\mathbb{I}_{(m+d)\times(m+d)})$ and setting $(\textbf{y},\textbf{x})^\top=L\textbf{z}$, where $LL^\top=C_{\xi;\eta}(\Delta)$. Consequently, $(\textbf{y},\textbf{x})^\top$ is a realization of the joint random vector 
\begin{equation*}
    \begin{pmatrix}
        \xi_i(\Delta) \\
        \eta_i(\Delta)
    \end{pmatrix},
\end{equation*}
where $\eta_i(\Delta)$ and $\xi_i(\Delta)$ are as in \eqref{eq:eta_i} and \eqref{eq:xi_i}, respectively.

When $d=m$, a natural alternative procedure (motivated by the standard approach mentioned in the beginning of Section \ref{sec2:LinearSDEs}) is to sample only a single $d$-dimensional realization $\textbf{z}$ from $\mathcal{N}_d(0_d,\mathbb{I}_{d\times d})$ and to set $\textbf{y}=\sqrt{\Delta}\textbf{z}$ and $\textbf{x}=L_\eta(\Delta) \textbf{z}$, where $L_\eta(\Delta) L_\eta(\Delta)^\top=C_\eta(\Delta)$. In this way, $\textbf{x}$ and $\textbf{y}$ are still sampled from the correct marginal distributions, that is $\textbf{x}$ and $\textbf{y}$ are realizations of $\eta_i(\Delta)\sim \mathcal{N}_d(0_d,C_\eta(\Delta))$ and $\xi_i(\Delta)\sim \mathcal{N}_d(0_d,C_\xi(\Delta))$, respectively, while both are generated from the same underlying random source. 

Unlike in Algorithm~\ref{alg:Exact_WPB}, however, this construction does not preserve the exact cross-covariance $C(\Delta)$ from \eqref{eq:covMatrix_C}. Indeed, letting $\textbf{Z}\sim \mathcal{N}_d(0_d,\mathbb{I}_{d\times d})$, 
\begin{equation*}
    \underline{C}(\Delta):=\textrm{Cov}(L_\eta(\Delta)\mathbf{Z},\sqrt{\Delta}\mathbf{Z})=\sqrt{\Delta} L_\eta(\Delta) \neq C(\Delta).
\end{equation*}
Nevertheless, in Lemmas \ref{lemma:d_m_1} and \ref{lemma:d_m_larger_1} below, we show that the difference of $\underline{C}(\Delta)$ and $C(\Delta)$ is of order $\Delta^3$ for $d=m=1$, and, under suitable assumptions on the diffusion matrix $\Sigma$, of order $\Delta^2$ for $d=m>1$.

\paragraph{Alternative exact simulation algorithm.} 

The alternative simulation procedure described above can be slightly modified for direct comparability with Algorithm \ref{alg:Exact_WPB}. Specifically, given a Wiener path $W(t_i)$, $i=0,\ldots,n$ (e.g., generated by Algorithm \ref{alg:Exact_WPB}), one may set 
\begin{equation*}
    \mathbf{y}=W(t_i)-W(t_{i-1}), \ \mathbf{z}=\mathbf{y}/\sqrt{\Delta}, \ \text{and} \ \mathbf{x}=L_\eta(\Delta)\mathbf{z}.
\end{equation*}
In this way, $\textbf{x}$ and $\textbf{y}$ are realizations of the marginal random vectors $\eta_i(\Delta)\sim \mathcal{N}_d(0_d,C_\eta(\Delta))$ and $\xi_i(\Delta)\sim\mathcal{N}_d(0_d,C_\xi(\Delta))$, respectively, while both are based on the same underlying random source. However, while Algorithm \ref{alg:Exact_WPB} generates $X(t_i)$ and $W(t_i)$, $i=0,\ldots,n$, with the correct cross-covariance structure, the alternative procedure, as detailed in Algorithm \ref{alg:Exact_alternative}, generates a path $\underline{X}(t_i)$, $i=0,\ldots,n$, which is linked to $W(t_i)$, $i=0,\ldots,n$, only via an approximate cross-covariance structure. 

Nevertheless, in Theorems \ref{prop:d_m_1} and \ref{prop:d_m_larger_1} below, we show that for any $n$ and $i=0,\ldots,n$, 
\begin{equation}\label{eq:Cerror}
    C_{\textrm{error}}(\Delta):=\textrm{Cov}(\underline{X}(t_i),W(t_i))-\textrm{Cov}(X(t_i),W(t_i))
\end{equation}
is of order $\Delta^2$ for $d=m=1$, and, under suitable assumptions on $A$, $\Sigma$, and $\Delta$, of order $\Delta$ for $d=m>1$. This means that, in the alternative construction, the cross-covariance structure between the simulated Wiener path and the corresponding solution path approximates the exact one with an error controlled by the time step $\Delta$, which tends to zero as $\Delta\to 0$. 

\begin{algorithm}[t]
	\caption{\textbf{Alternative} exact simulation (\textbf{Alternative-exact}) 
		\ \\ \textbf{Input:} Initial value $X_0$, time horizon $T$, 
  step-size $\Delta$, exponential matrix~$e^{A \Delta}$, covariance matrix $C_{\eta}(\Delta)$ and Wiener path $W(t_i)$, $i=0,\ldots,n$. 
  \ \\
		\textbf{Output:} Path of the linear SDE \eqref{eq:LinearSDE}.   
        }
	\label{alg:Exact_alternative}
    \begin{spacing}{1.1}
	\begin{algorithmic}[1]
        \State Set $n=T/\Delta$.
		\State Set $\underline{X}(t_0)=X_0$. 
        \State Determine $L_\eta(\Delta)$ such that $L_\eta(\Delta) L_\eta(\Delta)^\top=C_{\eta}(\Delta)$.
		\For{$i = 1:n$}
        \State Set $\mathbf{y}=W(t_i)-W(t_{i-1})$
        \State Set $\mathbf{z}=\mathbf{y}/\sqrt{\Delta}$
        \State Set $\mathbf{x}=L_\eta(\Delta)\mathbf{z}$.
		\State Set $\underline{X}(t_i)=e^{A\Delta}\underline{X}(t_{i-1})+\mathbf{x}$.  
		\EndFor
		\State Return $\underline{X}(t_i)$, 
        $i=0,\ldots,n$. 
	\end{algorithmic}
    \end{spacing}
\end{algorithm}

\begin{remark}
    The path $\underline{X}(t_i)$, $i=0,\ldots,n$, generated by Algorithm \ref{alg:Exact_alternative} is an exact realization of SDE \eqref{eq:LinearSDE}, in the sense that it preserves the joint distribution of the solution at the grid points. However, it is related to the given Wiener path $W(t_i)$, $i=0,\ldots,n$, only through an approximate cross-covariance structure. Consequently, although $\underline{X}(t_i)$ is generated according to the correct distribution, it generally differs from the realization $X(t_i)$ associated with the given Wiener path.
\end{remark}

\subsubsection{The case $d=m=1$}\label{sec:d_m_1}

In this subsection, we consider the $(d=1)$-dimensional SDE \eqref{eq:LinearSDE}, where $A=-a \neq 0$, $\Sigma=\sigma>0$ and $(W(t))_{t\in[0,T]}$ is an $(m=1)$-dimensional Wiener process. 

The following result on the local cross-covariance error holds.

\begin{lemma}\label{lemma:d_m_1}
    Let $C(\Delta)$ and $C_\eta(\Delta)$ be as in \eqref{eq:C_xi_eta} and define $\underline{C}(\Delta):=\sqrt{\Delta}\sqrt{C_\eta(\Delta)}$. Then, it holds that
    \begin{equation*}
        \underline{C}(\Delta)-C(\Delta)=\mathcal{O}(\Delta^3).
    \end{equation*}
\end{lemma}

\begin{proof}
    Define $r_k(\Delta):=\mathcal{O}(\Delta^k)$.

    First, note that $C(\Delta)$ can be written as
    \begin{equation*}
        C(\Delta)=\sigma \int_0^\Delta e^{-a(\Delta-s)} \ ds, 
    \end{equation*}
    and that
    \begin{equation*}
        e^{-a(\Delta-s)}=1-a(\Delta-s)+\frac{a^2(\Delta-s)^2}{2}+r_3(\Delta). 
    \end{equation*}
    This implies that
    \begin{align}\label{eq:C_expand}
        \nonumber C(\Delta)&=\sigma \int_0^\Delta \left( 1-a(\Delta-s) + \frac{a^2(\Delta-s)^2}{2} + r_3(\Delta) \right) \ ds \\
        &= \sigma \left( \Delta - \frac{a}{2}\Delta^2 +\frac{a^2}{6} \Delta^3 +r_4(\Delta) \right).
    \end{align}

    Similar to before, $C_\eta(\Delta)$ can be written as
    \begin{equation*}
        C_\eta(\Delta)=\sigma^2 \int_0^\Delta e^{-2a(\Delta-s)} \ ds.
    \end{equation*}
    Using that
    \begin{equation*}
        e^{-2a(\Delta-s)}=1-2a(\Delta-s)+2a^2(\Delta-s)^2+r_3(\Delta), 
    \end{equation*}
    implies that
    \begin{align}\label{eq:C_eta_expand1}
        \nonumber C_\eta(\Delta)&=\sigma^2 \int_0^\Delta \left( 1-2a(\Delta-s) + 2a^2(\Delta-s)^2 + r_3(\Delta) \right) \ ds \\
        &= \sigma^2 \left( \Delta - a\Delta^2 +\frac{2a^2}{3} \Delta^3 +r_4(\Delta) \right).
    \end{align}

    Now, since $C_\eta(\Delta)=\mathcal{O}(\Delta)$, its square root should be of order $\sqrt{\Delta}$. Therefore, we make the following ansatz:
    \begin{equation*}
        \sqrt{C_\eta(\Delta)}=\sigma \sqrt{\Delta} + x \Delta^{3/2} + y \Delta^{5/2} + r_{7/2}(\Delta),
    \end{equation*}
    where $x$ and $y$ are unknown constants to be determined. Squaring the ansatz yields
    \begin{align}\label{eq:C_eta_expand2}
        C_\eta(\Delta) 
        &=\sigma^2 \Delta + 2\sigma x \Delta^2 + (x^2+2\sigma y)\Delta^3+r_4(\Delta).
    \end{align}
    Comparing the coefficients in \eqref{eq:C_eta_expand1} and \eqref{eq:C_eta_expand2} yields 
    \begin{equation*}
        x=-\frac{a\sigma}{2}, \qquad y=\frac{5a^2\sigma}{24}.
    \end{equation*}
    The final expansion is thus given by
    \begin{equation*}
        \sqrt{C_\eta(\Delta)}=\sigma \sqrt{\Delta}-\frac{a\sigma}{2}\Delta^{3/2}+\frac{5a^2\sigma}{24}\Delta^{5/2}+r_{7/2}(\Delta),
    \end{equation*}
    and multiplying it by $\sqrt{\Delta}$ gives
    \begin{equation}\label{eq:C_underline_expand}
        \underline{C}(\Delta)=\sqrt{\Delta}\sqrt{C_\eta(\Delta)}=\sigma \Delta - \frac{a\sigma}{2} \Delta^2 +\frac{5 a^2 \sigma}{24}\Delta^3 + r_4(\Delta).
    \end{equation}

    Finally, taking the difference of \eqref{eq:C_expand} and \eqref{eq:C_underline_expand}, yields the result
    \begin{equation*}
        \underline C(\Delta)-C(\Delta)=\frac{a^2\sigma}{24}\Delta^3+r_4(\Delta)=\mathcal{O}(\Delta^3).
    \end{equation*}
\end{proof}

\begin{remark}
    As an alternative to the ansatz approach used in the proof of Lemma \ref{lemma:d_m_1}, one could derive \eqref{eq:C_underline_expand} by applying a suitable Taylor expansion of the square root function.
\end{remark}

Lemma \ref{lemma:d_m_1} shows that each time step introduces a local error of order $\Delta^3$ in the cross-covariance between the stochastic integral $\eta_i(\Delta)$ and the Wiener increment $\xi_i(\Delta)$. In the following theorem, we show that, after back-iterating the process to obtain $\underline{X}(t_i)$ at any grid point $t_i \in [0,T]$, this local error is amplified by the factor 
\begin{equation*}
    \sum\limits_{k=1}^i e^{-a(i-k)\Delta} = \frac{1-e^{- a t_i}}{1-e^{-a\Delta}}=\mathcal{O}(\Delta^{-1}),
\end{equation*}
resulting in a global error in the cross-covariance of order $\Delta^2$.

\begin{theorem}\label{prop:d_m_1}
    Let $(X(t_i),W(t_i))^\top$ and $(\underline{X}(t_i),W(t_i))^\top$, $i=0,\ldots,n$, denote the random vectors at the grid points generated by Algorithm~\ref{alg:Exact_WPB} and Algorithm~\ref{alg:Exact_alternative}, respectively, and let $C_{\textrm{error}}(\Delta)$ be as in \eqref{eq:Cerror}. Then, for any grid point $t_i \in [0,T]$, it holds that
    \begin{equation*}
        C_{\textrm{error}}(\Delta) = \mathcal{O}(\Delta^2).
    \end{equation*}
\end{theorem}

\begin{proof}
    First, note that for any grid point $t_i \in [0,T]$ one has
    \begin{equation*}
        X(t_i)=e^{-a t_i}X_0 + \sum\limits_{k=1}^i e^{-a(i-k)\Delta} \eta_k(\Delta), 
    \end{equation*}
    where
    \begin{equation*}
        \eta_k(\Delta)=\sigma \int\limits_{t_{k-1}}^{t_k} e^{-a(t_k-s)} \ dW(s) \sim \mathcal{N}(0,C_\eta(\Delta)),
    \end{equation*}
    and
    \begin{equation*}
        W(t_i)=\sum\limits_{k=1}^{i} \xi_k(\Delta),
    \end{equation*}
    where
    \begin{equation*}
        \xi_k(\Delta)=\int\limits_{t_{k-1}}^{t_k} 1 \ dW(s) =W(t_k) -W(t_{k-1}) \sim \mathcal{N}(0,C_\xi(\Delta)).
    \end{equation*}
    
    Moreover, since $X_0$ is independent of the Wiener process,
    \begin{equation*}
        \textrm{Cov}(e^{-a t_i}X_0,W(t_i))=0.
    \end{equation*} 
    In addition, since different Wiener increments are independent, it holds that
    \begin{equation*}
        \textrm{Cov}(\eta_k(\Delta),\xi_l(\Delta))=0, \quad k\neq l.
    \end{equation*}
    Therefore, 
    \begin{equation*}
        \textrm{Cov}(X(t_i),W(t_i))=\sum\limits_{k=1}^{i} e^{-a(i-k)\Delta} \textrm{Cov}(\eta_k(\Delta),\xi_k(\Delta)) = C(\Delta) \sum\limits_{k=1}^i e^{-a(i-k)\Delta},
    \end{equation*}
    and
    \begin{equation*}
        \textrm{Cov}(\underline{X}(t_i),W(t_i))= \underline{C}(\Delta) \sum\limits_{k=1}^i e^{-a(i-k)\Delta},
    \end{equation*}
    respectively.

    The global cross-covariance error is thus given by
    \begin{align*}
    C_{\textrm{error}}(\Delta) &:=\textrm{Cov}(\underline{X}(t_i),W(t_i))-\textrm{Cov}(X(t_i),W(t_i)) \\ &= \left( \underline{C}(\Delta) - C(\Delta) \right) \sum\limits_{k=1}^i e^{-a(i-k)\Delta} 
    =\left( \underline{C}(\Delta) - C(\Delta) \right)\frac{1-e^{-a t_i}}{1-e^{-a\Delta}}.
    \end{align*}
    Since
    \begin{equation*}
        1-e^{-a\Delta} = a\Delta + \mathcal{O}(\Delta^2),
    \end{equation*}
    we have
    \begin{equation*}
        \frac{1-e^{-a t_i}}{1-e^{-a\Delta}} = \mathcal{O}(\Delta^{-1}).
    \end{equation*}
    Together with Lemma \ref{lemma:d_m_1}, which yields
    \begin{equation*}
        \underline{C}(\Delta) - C(\Delta) = \mathcal{O}(\Delta^3),
    \end{equation*}
   this implies
    \begin{equation*}
        C_{\textrm{error}}(\Delta)=\mathcal{O}(\Delta^{-1})\mathcal{O}(\Delta^3)=\mathcal{O}(\Delta^2).
    \end{equation*}
\end{proof}

\subsubsection{The case $d=m>1$}\label{sec:d_m_larger1}

Let $d>1$. In this subsection, we consider the $d$-dimensional SDE \eqref{eq:LinearSDE}, where $A,\Sigma \in \mathbb{R}^{d\times d}$ and $(W(t))_{t\in[0,T]}$ is a $d$-dimensional Wiener process. 

The following result on the local cross-covariance error holds.

\begin{lemma}\label{lemma:d_m_larger_1}
    Assume that the diffusion matrix $\Sigma \in \mathbb{R}^{d\times d}$ is a lower triangular matrix with positive diagonal entries.  Moreover, let $C(\Delta)$ and $C_\eta(\Delta)$ be as in \eqref{eq:C_xi_eta} and define $\underline{C}(\Delta):=\sqrt{\Delta}L_\eta(\Delta)$, where $L_\eta(\Delta)L_\eta(\Delta)^\top=C_\eta(\Delta)$. Then, it holds that
    \begin{equation*}
        \underline{C}(\Delta)-C(\Delta)=\mathcal{O}(\Delta^2).
    \end{equation*}
\end{lemma}

\begin{proof}
    Define $R_k(\Delta):=\mathcal{O}(\Delta^k)$. 

    Using that
    \begin{equation*}
        e^{A(\Delta-s)}=\mathbb{I}_{d\times d}+A(\Delta-s)+R_2(\Delta),
    \end{equation*}
    yields
    \begin{align}\label{eq:C_expand_mult}
        \nonumber C(\Delta)&=\int_0^\Delta e^{A(\Delta-s)} \Sigma \ ds \\
        \nonumber &= \int\limits_0^\Delta \left( \mathbb{I}_{d\times d}+A(\Delta-s)+R_2(\Delta) \right)\Sigma \ ds \\
        &=\Delta \Sigma +\frac{\Delta^2}{2}A\Sigma + R_3(\Delta).
    \end{align}

    Moreover, we have that
    \begin{equation*}
        e^{A(\Delta-s)}\Sigma\Sigma^\top \left( e^{A(\Delta-s)} \right)^\top = \Sigma\Sigma^\top + (\Delta-s)\left( A \Sigma\Sigma^\top + \Sigma\Sigma^\top A^\top \right) + R_2(\Delta).
    \end{equation*}
    This implies that
    \begin{samepage}
        \begin{align}\label{eq:C_eta_expand1_mult}
        \nonumber C_\eta(\Delta)&= \int\limits_0^\Delta e^{A(\Delta-s)}\Sigma\Sigma^\top \left( e^{A(\Delta-s)} \right)^\top ds \\
        \nonumber &=\int\limits_0^\Delta \left( \Sigma\Sigma^\top + (\Delta-s)\left( A \Sigma\Sigma^\top + \Sigma\Sigma^\top A^\top \right) + R_2(\Delta) \right) ds \\
        &= \Delta \Sigma\Sigma^\top + \frac{\Delta^2}{2}\left( A \Sigma\Sigma^\top + \Sigma\Sigma^\top A^\top \right) + R_3(\Delta).
    \end{align}
    \end{samepage}    

    From \eqref{eq:C_eta_expand1_mult} we get that  $L_\eta(\Delta)L_\eta(\Delta)^\top=C_\eta(\Delta)=\mathcal{O}(\Delta)$.
    Moreover, by assumption $\Sigma$ is a lower triangular matrix with positive diagonal entries, and therefore $\Sigma$ is the unique Cholesky factor of~$\Sigma\Sigma^\top$.  Hence, we seek an expansion of $L_\eta(\Delta)$ of the form:
    \begin{equation}\label{eq:L_eta_expand}
        L_\eta(\Delta)= \sqrt{\Delta}\Sigma + \Delta^{3/2} X + R_{5/2}(\Delta),
    \end{equation}
    where $X$ is an unknown constant matrix. This implies that
    \begin{equation}\label{eq:C_eta_expand2_mult}
        C_\eta(\Delta)=L_\eta(\Delta)L_\eta(\Delta)^\top  = \Delta \Sigma\Sigma^\top + \Delta^2 \left( \Sigma X^\top + X \Sigma^\top \right) + R_3(\Delta).
    \end{equation}
    
    Comparing the coefficients in \eqref{eq:C_eta_expand1_mult} and \eqref{eq:C_eta_expand2_mult} leads to the matrix equation 
    \begin{equation}\label{eq:matrix_eq}
        \Sigma X^\top + X \Sigma^\top =\frac{1}{2}\left( A \Sigma\Sigma^\top + \Sigma \Sigma^\top A^\top \right). 
    \end{equation}
    It is straightforward to verify that $X_0=\frac{1}{2}A\Sigma$ is a particular solution of \eqref{eq:matrix_eq}. If $X$ is any solution of \eqref{eq:matrix_eq}, then the difference $Y:=X-X_0$ satisfies the homogeneous matrix equation 
    \begin{equation}\label{eq:hom_mat_eq.}
        \Sigma Y^\top+Y \Sigma^\top=\mathbb{O}_{d \times d}.
    \end{equation}
    Conversely, adding any solution $Y$ of this homogeneous equation to $X_0$ yields a solution of \eqref{eq:matrix_eq}. Hence, every solution is of the form $X=\frac{1}{2}A\Sigma+Y$, where $Y$ satisfies \eqref{eq:hom_mat_eq.}.
    
    Using \eqref{eq:L_eta_expand}, the final expansion is thus given by
    \begin{equation*}
        L_\eta(\Delta)=\sqrt{\Delta}\Sigma + \Delta^{3/2}\left( \frac{1}{2}A\Sigma +Y \right) + R_{5/2}(\Delta),
    \end{equation*}
    and multiplying by $\sqrt{\Delta}$ gives
    \begin{equation}\label{eq:underlineC_expand_mult}
        \underline{C}(\Delta)=\sqrt{\Delta}L_\eta(\Delta)=\Delta \Sigma + \Delta^2 \left( \frac{1}{2}A\Sigma +Y \right) + R_3(\Delta).
    \end{equation} 
    
    Consequently, taking the difference of \eqref{eq:C_expand_mult} and \eqref{eq:underlineC_expand_mult}, we obtain
    \begin{equation*}
        \underline{C}(\Delta)-C(\Delta)=\Delta^2 Y +R_3(\Delta) =\mathcal{O}(\Delta^2),
    \end{equation*}
    which proves the statement.
\end{proof}

\begin{remark}
    Note that, unlike in the scalar case $d=m=1$, the $\mathcal{O}(\Delta^2)$-term in the local cross-covariance error does not cancel in general, due to the freedom in choosing the matrix $Y$ as part of the solution to the matrix equation \eqref{eq:matrix_eq}. While $X_0=\frac{1}{2}A\Sigma$ solves the equation, it is not the unique solution. For example, let $\Sigma= \sigma \mathbb{I}_{2\times 2}$, $\sigma>0$ and 
    \begin{equation*}
        Y=\begin{bmatrix}
            0 & c \\
            -c & 0
        \end{bmatrix}, \quad c \in \mathbb{R}.
    \end{equation*}
    Then, $X=\frac{1}{2}A\Sigma+Y$ also solves the corresponding matrix equation \eqref{eq:matrix_eq}.
\end{remark}

Lemma \ref{lemma:d_m_larger_1} shows that, if $\Sigma$ is a lower triangular matrix with positive diagonal entries, each time step introduces a local error of order $\Delta^2$ in the cross-covariance between the stochastic integral $\eta_i(\Delta)$ and the Wiener increment $\xi_i(\Delta)$. Similar as in the scalar case, in the following theorem we show that, after back-iterating the process to obtain $\underline{X}(t_i)$ at any grid point $t_i \in [0,T]$, this local error is amplified by the factor 
\begin{equation*}
    \sum\limits_{k=1}^i e^{A(i-k)\Delta} = (\mathbb{I}_{d\times d}-e^{A t_i})(\mathbb{I}_{d\times d}-e^{A\Delta})^{-1}=\mathcal{O}(\Delta^{-1}).
\end{equation*}
Thus, requiring that $A$ and $\Delta$ are such that the matrix $\mathbb{I}_{d\times d}-e^{A\Delta}$ is invertible, the resulting global error in the cross-covariance is of order $\Delta$.

\begin{theorem}\label{prop:d_m_larger_1}
    Assume that the diffusion matrix $\Sigma \in \mathbb{R}^{d\times d}$ is a lower triangular matrix with positive diagonal entries, and that the drift matrix $A$ and the time step $\Delta>0$ are such that the matrix $\mathbb{I}_{d\times d}-e^{A\Delta}$ is invertible. Moreover, let $(X(t_i),W(t_i))^\top$ and $(\underline{X}(t_i),W(t_i))^\top$, $i=0,\ldots,n$, denote the random vectors at the grid points generated by Algorithm~\ref{alg:Exact_WPB} and Algorithm~\ref{alg:Exact_alternative}, respectively, and let $C_{\textrm{error}}(\Delta)$ be as in \eqref{eq:Cerror}. Then, for any grid point $t_i \in [0,T]$, it holds that
    \begin{equation*}
        C_{\textrm{error}}(\Delta) 
        = \mathcal{O}(\Delta).
    \end{equation*} 
\end{theorem}

\begin{proof}
    Similar as in the scalar case, for any grid point $t_i \in [0,T]$ one has
    \begin{equation*}
        X(t_i)=e^{A t_i}X_0 + \sum\limits_{k=1}^i e^{A(i-k)\Delta} \eta_k(\Delta), 
    \end{equation*}
    where
    \begin{equation*}
        \eta_k(\Delta)= \int\limits_{t_{k-1}}^{t_k} e^{A(t_k-s)} \Sigma \ dW(s) \sim \mathcal{N}_d(0_d,C_\eta(\Delta)),
    \end{equation*}
    and
    \begin{equation*}
        W(t_i)=\sum\limits_{k=1}^{i} \xi_k(\Delta),
    \end{equation*}
    where
    \begin{equation*}
        \xi_k(\Delta)=\int\limits_{t_{k-1}}^{t_k} \mathbb{I}_{d\times d} \ dW(s) =W(t_k) -W(t_{k-1}) \sim \mathcal{N}_d(0_d,C_\xi(\Delta)).
    \end{equation*}
    
    Since $X_0$ is independent of the Wiener process,
    \begin{equation*}
        \textrm{Cov}(e^{A t_i}X_0,W(t_i))=\mathbb{O}_{d\times d},
    \end{equation*} 
    and, since different Wiener increments are independent, 
    \begin{equation*}
        \textrm{Cov}(\eta_k(\Delta),\xi_l(\Delta))=\mathbb{O}_{d\times d}, \quad k\neq l.
    \end{equation*}
    Therefore, 
    \begin{equation*}
        \textrm{Cov}(X(t_i),W(t_i))=\sum\limits_{k=1}^{i} e^{A(i-k)\Delta} \textrm{Cov}(\eta_k(\Delta),\xi_k(\Delta)) = \sum\limits_{k=1}^i e^{A(i-k)\Delta}C(\Delta),
    \end{equation*}
    and
    \begin{equation*}
        \textrm{Cov}(\underline{X}(t_i),W(t_i))= \sum\limits_{k=1}^i e^{A(i-k)\Delta}\underline{C}(\Delta),
    \end{equation*}
    respectively.

    The global cross-covariance error is thus given by
    \begin{align*}
    C_{\textrm{error}}(\Delta) :=\textrm{Cov}(\underline{X}(t_i),W(t_i))-\textrm{Cov}(X(t_i),W(t_i)) = \sum\limits_{k=1}^i e^{A(i-k)\Delta} \left( \underline{C}(\Delta) - C(\Delta) \right).
    \end{align*}
    Now, since the matrix $\mathbb{I}_{d\times d}-e^{A\Delta}$ is invertible by assumption, one has
    \begin{align*}
    \sum\limits_{k=1}^i e^{A(i-k)\Delta} = (\mathbb{I}_{d\times d}-e^{A\Delta})^{-1} (\mathbb{I}_{d\times d}-e^{A t_i}).
    \end{align*}   
    Moreover, since 
    \begin{equation*}
        \mathbb{I}_{d\times d}-e^{A\Delta} = -A\Delta + \mathcal{O}(\Delta^2),
    \end{equation*}
    and, since $\mathbb{I}_{d\times d}-e^{A\Delta}$ is invertible, it follows that
    \begin{equation*}
      (\mathbb{I}_{d\times d}-e^{A\Delta})^{-1} (\mathbb{I}_{d\times d}-e^{A t_i}) = \mathcal{O}(\Delta^{-1}).
    \end{equation*}
    In addition, under the assumption that $\Sigma$ is a lower triangular matrix with positive diagonal entries, Lemma \ref{lemma:d_m_larger_1} further yields that
    \begin{equation*}
        \underline{C}(\Delta) - C(\Delta) = \mathcal{O}(\Delta^2).
    \end{equation*}
   Consequently,
    \begin{equation*}
        C_{\textrm{error}}(\Delta)=\mathcal{O}(\Delta^{-1})\mathcal{O}(\Delta^2)=\mathcal{O}(\Delta),
    \end{equation*}
    which proves the statement.
\end{proof}

\subsubsection{Illustration of the results on the Ornstein--Uhlenbeck process and simple stochastic harmonic oscillator} 

The theoretical results of Section \ref{sec:d_m_1} are illustrated on the OU process
\begin{equation*}
    dX(t)=-aX(t) dt + \sigma dW(t), \quad X(0)=X_0,
\end{equation*}
where we set $X_0=1$, $a=2$, and $\sigma=1/2$. The results of Section \ref{sec:d_m_larger1} are illustrated on the simple SHO introduced in Section \ref{sec3:simple_SHO}. For this example, we set $X_0=(0,0)^\top$, $\sigma_1=1$, and $\sigma_2=2$. Furthermore, we consider two choices for the remaining diffusion parameters. First, we set $\sigma_3=0$ and $\sigma_4=3/2$, such that $\Sigma$ is lower triangular with positive diagonal entries and therefore satisfies the assumption of Lemma \ref{lemma:d_m_larger_1} and Theorem \ref{prop:d_m_larger_1}. Second, we set $\sigma_3=3/2$ and $\sigma_4=0$, which violates this assumption. 

\begin{figure}[t]
	\begin{centering}
		\includegraphics[width=1.0\textwidth]{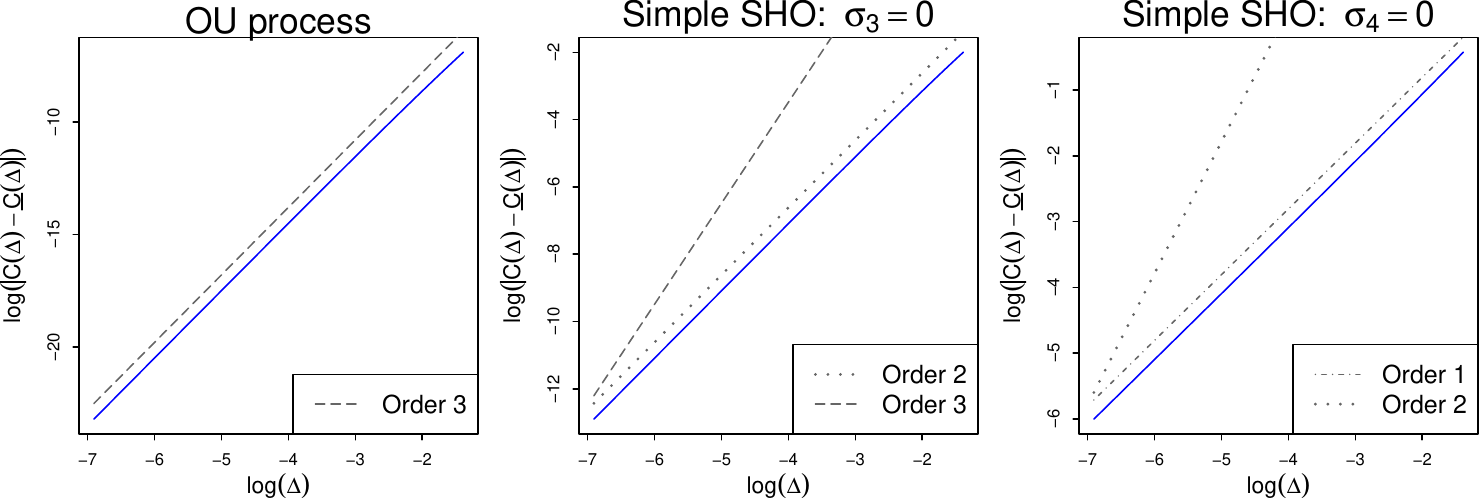}	
		\caption{Local cross-covariance error $|C(\Delta)-\underline{C}(\Delta)|$ as a function of the time step $\Delta$ in a log-log-plot. Left panel: OU process, showing the order $\Delta^3$ predicted by Lemma \ref{lemma:d_m_1}. Middle panel: simple SHO, where the assumption of Lemma \ref{lemma:d_m_larger_1} is satisfied, showing the predicted order $\Delta^2$. Right panel: simple SHO violating the assumption on $\Sigma$, where only order $\Delta$ is observed.}
		\label{fig:Local_cov_error}
	\end{centering}
\end{figure}

\begin{figure}[t]
	\begin{centering}
		\includegraphics[width=0.67\textwidth]{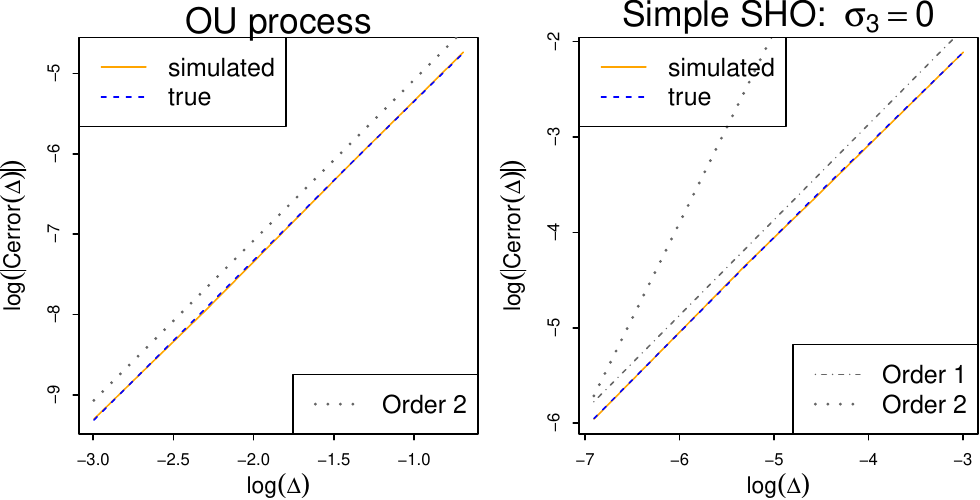}	
		\caption{Global cross-covariance error $C_\textrm{error}(\Delta)$ as a function of the time step $\Delta$ in a log-log-plot. The blue dashed lines correspond to the theoretical values, while the solid orange lines show the empirical estimates obtained from repeated path simulations using Algorithm \ref{alg:Exact_WPB} and Algorithm \ref{alg:Exact_alternative}. Lef panel: OU process, confirming the order $\Delta^2$ predicted by Theorem \ref{prop:d_m_1}. Right panel: simple SHO, where the assumptions of Theorem \ref{prop:d_m_larger_1} are satisfied, confirming the predicted order $\Delta$.}
		\label{fig:Global_cov_error}
	\end{centering}
\end{figure}

\begin{figure}[t]
	\begin{centering}
		\includegraphics[width=1.0\textwidth]{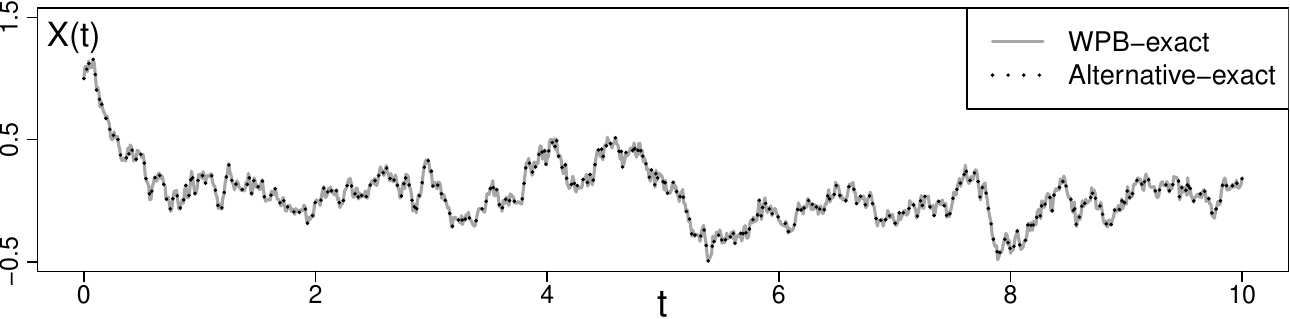}	
		\caption{Comparison of paths generated by Algorithm \ref{alg:Exact_WPB} and Algorithm \ref{alg:Exact_alternative} for the OU process. The solid gray line shows the path $X(t_i)$ generated with the exact cross-covariance structure, while the dotted black line shows the corresponding path $\underline{X}(t_i)$, generated using the alternative procedure based on the same Wiener path. 
        }
		\label{fig:OU_paths_cov_error}
	\end{centering}
\end{figure}

\begin{figure}[t]
	\begin{centering}
		\includegraphics[width=1.0\textwidth]{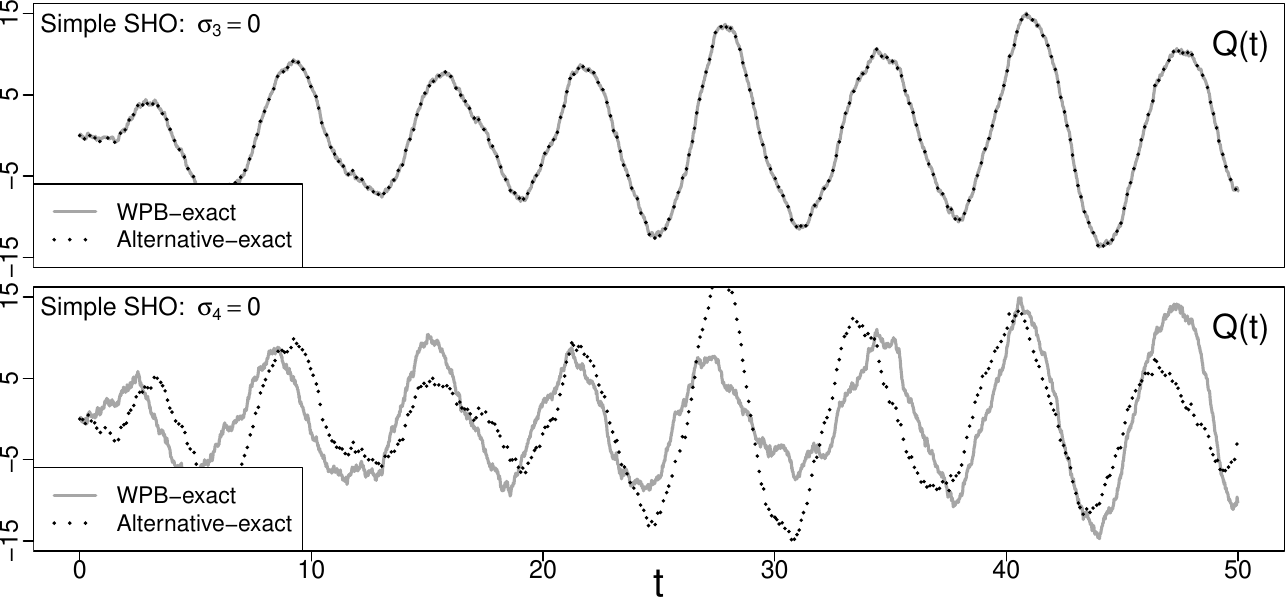}	
		\caption{Comparison of the $Q$-component of paths generated by Algorithm \ref{alg:Exact_WPB} and Algorithm \ref{alg:Exact_alternative} for the simple SHO. Top panel: The assumptions of Theorem \ref{prop:d_m_larger_1} are satisfied and the two paths closely agree. Bottom panel: The assumption on $\Sigma$ is violated and the two paths exhibit a visible deviation.}
		\label{fig:SHO_paths_cov_error}
	\end{centering}
\end{figure}

In Figure \ref{fig:Local_cov_error}, we compute the true local cross-covariance  $C(\Delta)$ and its approximation $\underline{C}(\Delta)$ as functions of the step size $\Delta$. The quantity $\log(\Delta)$ versus $\log(|C(\Delta)-\underline{C}(\Delta)|)$ is plotted (blue lines) to illustrate the order of the local cross-covariance error. For the OU process (left panel), we observe an order of $\Delta^3$, in agreement with Lemma \ref{lemma:d_m_1}. For the simple SHO with $\sigma_3=0$ and $\sigma_4=3/2$ (middle panel), the diffusion matrix $\Sigma$ satisfies the assumption of Lemma \ref{lemma:d_m_larger_1}, and the observed order is $\Delta^2$, as predicted. In contrast, when the assumption on $\Sigma$ is violated by choosing $\sigma_3=3/2$ and $\sigma_4=0$ (right panel), only an order of $\Delta$ is observed. As discussed above, this local order is insufficient to guarantee a global cross-covariance error that converges to zero as $\Delta$ tends to zero.

In Figure \ref{fig:Global_cov_error}, we consider the global cross-covariance error $C_\textrm{error}(\Delta)$ defined in \eqref{eq:Cerror} as a function of $\Delta$ and plot $\log(\Delta)$ versus $\log(|C_\textrm{error}(\Delta)|)$ (blue dashed lines). Specifically, we compare the theoretical values with estimates obtained from repeated path simulations using Algorithm \ref{alg:Exact_WPB} and Algorithm \ref{alg:Exact_alternative}, where the covariance is estimated empirically. The simulated estimates agree well with the theoretical values. 
For the OU process, we observe an order of $\Delta^2$, in agreement with Theorem \ref{prop:d_m_1}. For the simple SHO with $\sigma_3=0$ and $\sigma_4=3/2$, we observe the predicted order~$\Delta$ from Theorem \ref{prop:d_m_larger_1}. 

\begin{remark}
    Note that, for the simple SHO the additional assumption of Theorem \ref{prop:d_m_larger_1} that the matrix $\mathbb{I}_{2\times 2}-e^{A\Delta}$ is invertible is satisfied as long as $\Delta \neq 2\pi k$, $k\in \mathbb{N}$, since for its determinant it holds that $\textrm{det}(\mathbb{I}_{2\times 2}-e^{A\Delta})=2(1-\cos(\Delta))$. This does not impose a practical restriction, since the time steps considered in numerical simulations are typically much smaller than one.
\end{remark}

Finally, we compare sample paths generated by Algorithm \ref{alg:Exact_WPB} and Algorithm \ref{alg:Exact_alternative}. Figure \ref{fig:OU_paths_cov_error} shows a realization $X(t_i)$ (solid gray line) of the OU process together with the corresponding path $\underline{X}(t_i)$ generated by the alternative procedure (dotted black line). The two paths are visually indistinguishable and this observed agreement is in accordance with the convergence of the global cross-covariance error established in Theorem \ref{prop:d_m_1}. 

Figure \ref{fig:SHO_paths_cov_error} shows the corresponding comparison for the $Q$-component of the simple SHO. In the top panel, where $\Sigma$ satisfies the assumption of Theorem \ref{prop:d_m_larger_1}, the two paths again closely agree. In the bottom panel, where the assumption on $\Sigma$ is violated, the two paths exhibit a visible deviation. Hence, in this case, the alternative procedure does not provide a close approximation of the path associated with the given Wiener realization, although it still generates a path with the correct marginal distribution of the SDE solution.

\subsection{Semi-linear SDEs with additive noise for $d=m$}\label{sec6:SemiLinSDEs}

In the following, we extend Algorithm \ref{alg:Exact_alternative} to semi-linear SDEs \eqref{eq:semi-linearSDE} using the splitting schemes \eqref{eq:splitting}, and then illustrate this approach on the elliptic stochastic FHN model \eqref{eq:FHN_elliptic}. 

\paragraph{Alternative Strang simulation algorithm.}

Recall that the splitting schemes \eqref{eq:splitting} for SDE \eqref{eq:semi-linearSDE} contain the independent stochastic integrals $\eta_i(\Delta)$, $i=1,\ldots,n$, defined in \eqref{eq:eta_i}. Therefore, similar to before, if $d=m$ and a $d$-dimensional Wiener path $W(t_i)$, $i=0,\ldots,n$, is given, one may set
\begin{equation*}
    \mathbf{y}=W(t_i)-W(t_{i-1}), \ \mathbf{z}=\mathbf{y}/\sqrt{\Delta}, \ \text{and} \ \mathbf{x}=L_\eta(\Delta)\mathbf{z},
\end{equation*}
to obtain a realization $\mathbf{x}$ of $\eta_i(\Delta)$ that is based on the same underlying random source as the corresponding Wiener increment. Using this $\mathbf{x}$ to update the approximate SDE path via the splitting schemes \eqref{eq:splitting} yields Algorithm \ref{alg:Strang_alternative} (reported for the Strang scheme). This algorithm generates an approximate path $\underline{\widetilde{X}}^{\textrm{S}}(t_i)$, $i=0,\ldots,n$, of SDE \eqref{eq:semi-linearSDE}, which is linked to the given Wiener path $W(t_i)$, $i=0,\ldots,n$, via an approximate cross-covariance structure (cf. Section \ref{sec6:LinSDEs}).

\begin{algorithm}[t]
	\caption{\textbf{Alternative} Strang simulation (\textbf{Alternative-Strang}) 
		\ \\ \textbf{Input:} Initial value $X_0$, time horizon $T$, 
  step-size $\Delta$, exponential matrix~$e^{A \Delta}$, covariance matrix $C_{\eta}(\Delta)$, Wiener path $W(t_i)$, $i=0,\ldots,n$, and ODE solution function $f(\cdot;\cdot)$. 
  \ \\
		\textbf{Output:} Approximate path of the semi-linear SDE \eqref{eq:semi-linearSDE}.   
        }
	\label{alg:Strang_alternative}
    \begin{spacing}{1.1}
	\begin{algorithmic}[1]
        \State Set $n=T/\Delta$.
		\State Set $\underline{\widetilde{X}}^{\textrm{S}}(t_0)=X_0$. 
        \State Determine $L_\eta(\Delta)$ such that $L_\eta(\Delta) L_\eta(\Delta)^\top=C_{\eta}(\Delta)$.
		\For{$i = 1:n$}
        \State Set $\mathbf{y}=W(t_i)-W(t_{i-1})$
        \State Set $\mathbf{z}=\mathbf{y}/\sqrt{\Delta}$
        \State Set $\mathbf{x}=L_\eta(\Delta)\mathbf{z}$.
		\State Set $\underline{\widetilde{X}}^{\textrm{S}}(t_i)=f\Bigl(e^{A\Delta}f(\underline{\widetilde{X}}^{\textrm{S}}(t_{i-1});\Delta/2)+\mathbf{x};\Delta/2\Bigr)$.  
		\EndFor
		\State Return $\underline{\widetilde{X}}^{\textrm{S}}(t_i)$, 
        $i=0,\ldots,n$. 
	\end{algorithmic}
    \end{spacing}
\end{algorithm}

\paragraph{Illustration of the algorithm on the elliptic stochastic FHN model.}

This approach is now illustrated on the elliptic stochastic FHN model \eqref{eq:FHN_elliptic}, for which the diffusion matrix $\Sigma$ trivially satisfies the assumption of Theorem \ref{prop:d_m_larger_1}.  

\begin{remark}
    The additional assumption of Theorem \ref{prop:d_m_larger_1} that the matrix $\mathbb{I}_{2\times 2}-e^{A\Delta}$ (with drift matrix $A$ as in \eqref{eq:splitting_choice_FHN} and matrix exponential $e^{A\Delta}$ as in \eqref{eq:MatrixExp_FHN}) is invertible also holds, since 
    \begin{equation*}
        \textrm{det}(\mathbb{I}_{2\times 2}-e^{A\Delta})=\bigl(1-e^{-\Delta/2}\bigr)^2+2e^{-\Delta/2}\Bigl( 1-\cos(\sqrt{\kappa}\Delta/2) \Bigr)>0, \quad \text{for all } \Delta>0.
    \end{equation*}
\end{remark}

\begin{figure}[t]
	\begin{centering}
		\includegraphics[width=1.0\textwidth]{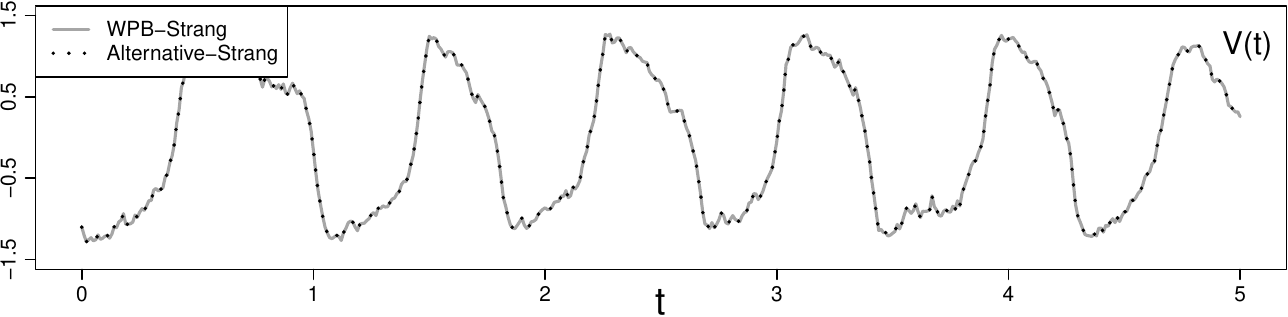}	
		\caption{Comparison of the $V$-component of paths generated by Algorithm \ref{alg:Strang_WPB} and Algorithm~\ref{alg:Strang_alternative} for the elliptic stochastic FHN model \eqref{eq:FHN_elliptic}. The two paths closely agree.} 
		\label{fig:FHN_paths_WPB_vs_Alternative}
	\end{centering}
\end{figure}

In the following example, the same input components and parameter values as in Section \ref{sec5.1:FHN_hypo} are used. Figure \ref{fig:FHN_paths_WPB_vs_Alternative} shows an approximate sample path $\widetilde{X}^{\textrm{S}}(t_i)$ (gray solid line; the same as in Figure \ref{fig:FHN_elliptic}) generated by Algorithm~\ref{alg:Strang_WPB} (\textit{WPB-Strang}) together with the corresponding path $\underline{\widetilde{X}}^{\textrm{S}}(t_i)$ (black dotted line) obtained by Algorithm~\ref{alg:Strang_alternative} (\textit{Alternative-Strang}). As expected from the analysis of the previous subsection, the two paths closely overlap and are visually indistinguishable. 

\section{Conclusion}
\label{sec:conclusion}

In this work, we developed Wiener-path-based simulation procedures for linear and semi-linear SDEs with additive noise that preserve the connection between the SDE path and the corresponding realization of its driving Wiener process.

Building upon the Wiener-path-based exact simulation approach for linear SDEs introduced in~\cite{Cruz2020}, we extended the construction to the simultaneous generation of coherent paths on different time grids. Based on Proposition \ref{prop:coarse_eta_xi}, the coarse-grid path can be constructed directly from the fine-grid random components while preserving the underlying Wiener-path realization. This provides a natural correspondence between different discretization levels and is particularly useful in applications where paths obtained with different time-step sizes need to be compared. 

We further proposed an alternative exact simulation procedure for linear SDEs in the case $d=m$, which requires sampling from a single $d$-dimensional normal distribution instead of jointly sampling the Wiener increment and the stochastic integral, which together form an $(m+d)$-dimensional normal random vector. While this construction preserves the correct marginal distributions of these quantities, it does not preserve their exact cross-covariance structure. Nevertheless, we showed that the deviation of the resulting cross-covariance between the solution and the driving Wiener process at the grid points from its exact value is of order $\Delta^2$ for $d=m=1$, as established in Theorem \ref{prop:d_m_1}, and of order $\Delta$ for $d=m>1$, as established in Theorem \ref{prop:d_m_larger_1}. Thus, the proposed alternative procedure provides a computationally more efficient construction by reducing the number of random numbers required, while retaining a controlled approximation of the coupling between the SDE path and its underlying Wiener path.

The proposed procedures were successfully illustrated on the simple stochastic harmonic oscillator (SHO), a well-known linear SDE arising in classical mechanics, and the numerical results confirmed the theoretical findings. In particular, when the diffusion matrix of the simple SHO satisfies the assumption of Theorem \ref{prop:d_m_larger_1}, the path generated by the alternative exact procedure closely agrees with the path generated by the Wiener-path-based exact procedure. In contrast, when this assumption is violated, the two paths differ noticeably, in agreement with the theoretical analysis. In the scalar case $d=m=1$, we additionally investigated the cross-covariance for the Ornstein--Uhlenbeck process, where the numerical results likewise confirmed the theoretically established convergence order.

We then considered the more general setting of semi-linear SDEs with additive noise, whose drift contains an additional non-linear function $N$ satisfying either a global Lipschitz condition or a one-sided Lipschitz condition together with a polynomial growth bound. Such SDEs arise in a wide range of relevant applications. Unlike in the linear case, their solutions are generally not available in closed form, and numerical approximation methods are therefore required. To this end, we employed structure-preserving splitting methods for semi-linear SDEs and extended the Wiener-path-based simulation procedures developed for linear SDEs to this setting. In particular, we incorporated the proposed Wiener-path-based constructions into splitting methods for both generating approximate SDE paths coherent with their underlying Wiener paths and simulating corresponding paths on different time grids. Moreover, the alternative  procedure for $d=m$ was adapted to splitting simulation. 

Finally, the resulting algorithms were successfully applied to the stochastic FitzHugh--Nagumo model and the stochastic Jansen--Rit neural mass model, two well-established models from neuroscience. The numerical experiments demonstrate the applicability of the Wiener-path-based simulation procedures to non-linear stochastic systems and highlight their potential for pathwise comparisons of different numerical methods across different time discretizations.

Overall, this work establishes a set of simulation procedures for generating paths of linear and semi-linear SDEs with a controlled connection to their underlying Wiener paths, enabling their comparison across different numerical methods and time grids.


\vspace{1.0cm}
\noindent\textbf{Declaration of generative AI use.} 

\medskip

\noindent During the preparation of this work the author used OpenAI ChatGPT to assist with language refinement, clarity and presentation. After using this tool, the author critically reviewed, verified and edited the content as needed. The author takes full responsibility for the content of this manuscript.




\end{document}